\documentclass[12pt]{article}

\usepackage{amsmath,amssymb,amsthm,amscd,a4wide}

\begin{document}

\newcommand{\End}{{\rm{End}\ts}}
\newcommand{\Hom}{{\rm{Hom}}}
\newcommand{\Mat}{{\rm{Mat}}}
\newcommand{\ch}{{\rm{ch}\ts}}
\newcommand{\chara}{{\rm{char}\ts}}
\newcommand{\diag}{ {\rm diag}}
\newcommand{\non}{\nonumber}
\newcommand{\wt}{\widetilde}
\newcommand{\wh}{\widehat}
\newcommand{\ot}{\otimes}
\newcommand{\la}{\lambda}
\newcommand{\La}{\Lambda}
\newcommand{\De}{\Delta}
\newcommand{\al}{\alpha}
\newcommand{\be}{\beta}
\newcommand{\ga}{\gamma}
\newcommand{\Ga}{\Gamma}
\newcommand{\ep}{\epsilon}
\newcommand{\ka}{\kappa}
\newcommand{\vk}{\varkappa}
\newcommand{\si}{\sigma}
\newcommand{\vp}{\varphi}
\newcommand{\de}{\delta}
\newcommand{\ze}{\zeta}
\newcommand{\om}{\omega}
\newcommand{\ee}{\epsilon^{}}
\newcommand{\su}{s^{}}
\newcommand{\hra}{\hookrightarrow}
\newcommand{\ve}{\varepsilon}
\newcommand{\ts}{\,}
\newcommand{\vac}{\mathbf{1}}
\newcommand{\di}{\partial}
\newcommand{\qin}{q^{-1}}
\newcommand{\tss}{\hspace{1pt}}
\newcommand{\Sr}{ {\rm S}}
\newcommand{\U}{ {\rm U}}
\newcommand{\BL}{ {\overline L}}
\newcommand{\BE}{ {\overline E}}
\newcommand{\BP}{ {\overline P}}
\newcommand{\AAb}{\mathbb{A}\tss}
\newcommand{\CC}{\mathbb{C}\tss}
\newcommand{\KK}{\mathbb{K}\tss}
\newcommand{\QQ}{\mathbb{Q}\tss}
\newcommand{\SSb}{\mathbb{S}\tss}
\newcommand{\ZZ}{\mathbb{Z}\tss}
\newcommand{\X}{ {\rm X}}
\newcommand{\Y}{ {\rm Y}}
\newcommand{\Z}{{\rm Z}}
\newcommand{\Ac}{\mathcal{A}}
\newcommand{\Lc}{\mathcal{L}}
\newcommand{\Mc}{\mathcal{M}}
\newcommand{\Pc}{\mathcal{P}}
\newcommand{\Qc}{\mathcal{Q}}
\newcommand{\Tc}{\mathcal{T}}
\newcommand{\Sc}{\mathcal{S}}
\newcommand{\Bc}{\mathcal{B}}
\newcommand{\Ec}{\mathcal{E}}
\newcommand{\Fc}{\mathcal{F}}
\newcommand{\Hc}{\mathcal{H}}
\newcommand{\Uc}{\mathcal{U}}
\newcommand{\Vc}{\mathcal{V}}
\newcommand{\Wc}{\mathcal{W}}
\newcommand{\Yc}{\mathcal{Y}}
\newcommand{\Ar}{{\rm A}}
\newcommand{\Br}{{\rm B}}
\newcommand{\Ir}{{\rm I}}
\newcommand{\Fr}{{\rm F}}
\newcommand{\Jr}{{\rm J}}
\newcommand{\Or}{{\rm O}}
\newcommand{\GL}{{\rm GL}}
\newcommand{\Spr}{{\rm Sp}}
\newcommand{\Rr}{{\rm R}}
\newcommand{\Zr}{{\rm Z}}
\newcommand{\gl}{\mathfrak{gl}}
\newcommand{\middd}{{\rm mid}}
\newcommand{\ev}{{\rm ev}}
\newcommand{\Pf}{{\rm Pf}}
\newcommand{\Norm}{{\rm Norm\tss}}
\newcommand{\oa}{\mathfrak{o}}
\newcommand{\spa}{\mathfrak{sp}}
\newcommand{\osp}{\mathfrak{osp}}
\newcommand{\g}{\mathfrak{g}}
\newcommand{\h}{\mathfrak h}
\newcommand{\n}{\mathfrak n}
\newcommand{\z}{\mathfrak{z}}
\newcommand{\Zgot}{\mathfrak{Z}}
\newcommand{\p}{\mathfrak{p}}
\newcommand{\sll}{\mathfrak{sl}}
\newcommand{\agot}{\mathfrak{a}}
\newcommand{\qdet}{ {\rm qdet}\ts}
\newcommand{\Ber}{ {\rm Ber}\ts}
\newcommand{\HC}{ {\mathcal HC}}
\newcommand{\cdet}{ {\rm cdet}}
\newcommand{\tr}{ {\rm tr}}
\newcommand{\gr}{ {\rm gr}}
\newcommand{\str}{ {\rm str}}
\newcommand{\loc}{{\rm loc}}
\newcommand{\Gr}{{\rm G}}
\newcommand{\sgn}{ {\rm sgn}\ts}
\newcommand{\ba}{\bar{a}}
\newcommand{\bb}{\bar{b}}
\newcommand{\bi}{\bar{\imath}}
\newcommand{\bj}{\bar{\jmath}}
\newcommand{\bk}{\bar{k}}
\newcommand{\bl}{\bar{l}}
\newcommand{\hb}{\mathbf{h}}
\newcommand{\Sym}{\mathfrak S}
\newcommand{\fand}{\quad\text{and}\quad}
\newcommand{\Fand}{\qquad\text{and}\qquad}
\newcommand{\For}{\qquad\text{or}\qquad}
\newcommand{\OR}{\qquad\text{or}\qquad}

\renewcommand{\theequation}{\arabic{section}.\arabic{equation}}

\newtheorem{thm}{Theorem}[section]
\newtheorem{lem}[thm]{Lemma}
\newtheorem{prop}[thm]{Proposition}
\newtheorem{cor}[thm]{Corollary}
\newtheorem{conj}[thm]{Conjecture}
\newtheorem*{mthm}{Main Theorem}
\newtheorem*{mthma}{Theorem A}
\newtheorem*{mthmb}{Theorem B}

\theoremstyle{definition}
\newtheorem{defin}[thm]{Definition}

\theoremstyle{remark}
\newtheorem{remark}[thm]{Remark}
\newtheorem{example}[thm]{Example}

\newcommand{\bth}{\begin{thm}}
\renewcommand{\eth}{\end{thm}}
\newcommand{\bpr}{\begin{prop}}
\newcommand{\epr}{\end{prop}}
\newcommand{\ble}{\begin{lem}}
\newcommand{\ele}{\end{lem}}
\newcommand{\bco}{\begin{cor}}
\newcommand{\eco}{\end{cor}}
\newcommand{\bde}{\begin{defin}}
\newcommand{\ede}{\end{defin}}
\newcommand{\bex}{\begin{example}}
\newcommand{\eex}{\end{example}}
\newcommand{\bre}{\begin{remark}}
\newcommand{\ere}{\end{remark}}
\newcommand{\bcj}{\begin{conj}}
\newcommand{\ecj}{\end{conj}}

\newcommand{\bal}{\begin{aligned}}
\newcommand{\eal}{\end{aligned}}
\newcommand{\beq}{\begin{equation}}
\newcommand{\eeq}{\end{equation}}
\newcommand{\ben}{\begin{equation*}}
\newcommand{\een}{\end{equation*}}

\newcommand{\bpf}{\begin{proof}}
\newcommand{\epf}{\end{proof}}

\def\beql#1{\begin{equation}\label{#1}}

\title{\Large\bf Yangian characters and classical $\Wc$-algebras}

\author{{A. I. Molev\quad and\quad
E. E. Mukhin}}

\date{} % Start July 2012
\maketitle

\vspace{30 mm}

\begin{abstract}
The Yangian characters (or $q$-characters) are known to be closely related
to the classical $\Wc$-algebras and to the centers of the affine
vertex algebras at the critical level.
We make this relationship more explicit by producing
families of generators of the $\Wc$-algebras from the
characters of the Kirillov--Reshetikhin modules associated
with multiples of the first
fundamental weight in types $B$ and $D$
and of the fundamental modules in type $C$.
We also give an independent derivation of the character formulas
for these representations
in the context of the $RTT$ presentation
of the Yangians.
In all cases the generators of the $\Wc$-algebras
correspond to the recently constructed elements of the Feigin--Frenkel
centers via an affine version of the Harish-Chandra isomorphism.

\end{abstract}

%%%\vspace{5 mm}
%%%
%%%{\it Key words:}
%%%

\vspace{35 mm}

\noindent
School of Mathematics and Statistics\newline
University of Sydney,
NSW 2006, Australia\newline
alexander.molev@sydney.edu.au

\vspace{7 mm}

\noindent
Department of Mathematical Sciences\newline
Indiana University -- Purdue University Indianapolis\newline
402 North Blackford St, Indianapolis, IN 46202-3216, USA\newline
mukhin@math.iupui.edu

\newpage

\tableofcontents

\newpage

\section{Introduction}
\label{sec:int}
\setcounter{equation}{0}

{\bf 1.1.}\quad
Let $\g$ be a simple Lie algebra over $\CC$.
Choose a Cartan subalgebra $\h$ of $\g$ and
a triangular decomposition $\g=\n_-\oplus\h\oplus \n_+$.
Recall that the {\it Harish-Chandra homomorphism\/}
\beql{hchcl}
\U(\g)^{\h}\to \U(\h)
\eeq
is the projection of the $\h$-centralizer $\U(\g)^{\h}$ in the universal enveloping
algebra to $\U(\h)$ whose kernel is the two-sided ideal
$\U(\g)^{\h}\cap \U(\g)\n_+$. The restriction of the homomorphism \eqref{hchcl}
to the center $\Z(\g)$ of $\U(\g)$ yields an isomorphism
\beql{hchicl}
\Z(\g)\to \U(\h)^W
\eeq
called the {\it Harish-Chandra isomorphism\/}, where $\U(\h)^W$ denotes the subalgebra
of invariants in $\U(\h)$ with respect to an action of the Weyl group $W$ of $\g$;
see e.g. \cite[Sec.~7.4]{d:ae}.

In this paper we will be concerned with an affine version of the
isomorphism \eqref{hchicl}.
Consider
the
affine Kac--Moody algebra $\wh\g$ which is the central extension
\ben
\wh\g=\g\tss[t,t^{-1}]\oplus\CC K,
\een
where $\g[t,t^{-1}]$ is the Lie algebra of Laurent
polynomials in $t$ with coefficients in $\g$.
We have a natural analogue of the homomorphism \eqref{hchcl},
\beql{hchaff}
\U\big(t^{-1}\g[t^{-1}]\big)^{\h}\to \U\big(t^{-1}\h[t^{-1}]\big).
\eeq
The
vacuum module $V_{-h^{\vee}}(\g)$
at the critical level over $\wh\g$ is defined as
the quotient of the
universal enveloping algebra $\U(\wh\g)$ by the left
ideal generated by $\g[t]$ and $K+h^{\vee}$, where $h^{\vee}$
denotes the dual Coxeter number for $\g$.
The vacuum module $V_{-h^{\vee}}(\g)$ possesses a vertex algebra structure;
see e.g. \cite[Ch.~2]{f:lc}. The {\it center\/} of this vertex algebra
is defined by
\ben
\z(\wh\g)
=\{S\in V_{-h^{\vee}}(\g)\ |\ \g[t]\ts S=0\},
\een
its elements are called {\it Segal--Sugawara vectors\/}.
The center is a commutative
associative algebra which can be regarded as
a commutative subalgebra of $\U\big(t^{-1}\g[t^{-1}]\big)^{\h}$.
By the results of Feigin and Frenkel~\cite{ff:ak},
$\z(\wh\g)$ is an algebra of polynomials in infinitely many variables
and the restriction
of the homomorphism \eqref{hchaff} to
the subalgebra $\z(\wh\g)$ yields an isomorphism
\beql{hchiaff}
\z(\wh\g)\to \Wc({}^L\g),
\eeq
where $\Wc({}^L\g)$ is the {\it classical $\Wc$-algebra\/} associated with the
Langlands dual Lie algebra ${}^L\g$; see \cite{f:lc} for
a detailed exposition of these results.
The $\Wc$-algebra $\Wc({}^L\g)$ can be defined as a subalgebra
of $\U\big(t^{-1}\h[t^{-1}]\big)$ which consists of the elements
annihilated by the {\it screening operators\/}; see Sec.~\ref{sec:cw} below.

Recently, explicit generators of the Feigin--Frenkel center $\z(\wh\g)$
were constructed for the Lie algebras $\g$ of all classical types $A$, $B$, $C$ and $D$;
see \cite{cm:ho}, \cite{ct:qs} and \cite{m:ff}.
Our aim in this paper is to describe the Harish-Chandra images of these generators
in types $B$, $C$ and $D$. The corresponding results in type $A$ are given in
\cite{cm:ho}; we also reproduce them below in a slightly different form
(as in \cite{cm:ho}, we work with the reductive Lie algebra $\gl_N$ rather than the simple
Lie algebra $\sll_N$ of type $A$).
The images of the generators of $\z(\wh\g)$
under the isomorphism \eqref{hchiaff} turn out to be elements
of the $\Wc$-algebra $\Wc({}^L\g)$ written in terms of noncommutative
analogues of the complete and elementary
symmetric functions.

In more detail, for any $X\in\g$ and $r\in\ZZ$ introduce
the corresponding elements of the loop algebra $\g[t,t^{-1}]$
by $X[r]=X\tss t^r$.
The extended Lie algebra $\wh\g\oplus\CC\tau$ with $\tau=-d/dt$ is defined by
the commutation relations
\beql{taur}
\big[\tau,X[r]\tss\big]=-r\ts X[r-1],\qquad
\big[\tau,K\big]=0.
\eeq
Consider the natural extension of \eqref{hchiaff}
to the isomorphism
\beql{hchiaffe}
\chi:\z(\wh\g)\ot\CC[\tau]\to \Wc({}^L\g)\ot\CC[\tau],
\eeq
which is identical on $\CC[\tau]$; see Sec.~\ref{sec:gw} for the definition of $\chi$.

\medskip
\noindent
{\bf 1.2.}\quad
First let $\g=\gl_N$ be the general linear Lie algebra
with the standard basis elements $E_{ij}$, $1\leqslant i,j\leqslant N$.
For each $a\in\{1,\dots,m\}$
introduce the element $E[r]_a$ of the algebra
\beql{tenprka}
\underbrace{\End\CC^{N}\ot\dots\ot\End\CC^{N}}_m{}\ot\U
\eeq
by
\beql{matnota}
E[r]_a=\sum_{i,j=1}^{N}
1^{\ot(a-1)}\ot e_{ij}\ot 1^{\ot(m-a)}\ot E_{ij}[r],
\eeq
where the $e_{ij}$ are the standard matrix units and $\U$ stands
for the universal enveloping algebra of
$\wh\gl_N\oplus\CC\tau$.
Let $H^{(m)}$ and $A^{(m)}$ denote the respective images of the
symmetrizer and anti-symmetrizer in the group algebra
for the symmetric group $\Sym_m$ under
its natural action on $(\CC^{N})^{\ot m}$; see \eqref{symasym}.
We will identify $H^{(m)}$ and $A^{(m)}$ with the elements
$H^{(m)}\ot 1$ and $A^{(m)}\ot 1$ of the algebra \eqref{tenprka}.
Define the elements
$\phi^{}_{m\tss a},\psi^{}_{m\tss a}\in\U\big(t^{-1}\gl_N[t^{-1}]\big)$
by the expansions
\begin{align}\label{deftra}
\tr\ts A^{(m)} \big(\tau+E[-1]_1\big)\dots \big(\tau+E[-1]_m\big)
&=\phi^{}_{m\tss0}\ts\tau^m+\phi^{}_{m\tss1}\ts\tau^{m-1}
+\dots+\phi^{}_{m\tss m},\\[0.7em]
\label{deftrh}
\tr\ts H^{(m)} \big(\tau+E[-1]_1\big)\dots \big(\tau+E[-1]_m\big)
&=\psi^{}_{m\tss0}\ts\tau^m+\psi^{}_{m\tss1}\ts\tau^{m-1}
+\dots+\psi^{}_{m\tss m},
\end{align}
where the traces are taken
over all $m$ copies of $\End\CC^{N}$.
The results of \cite{cm:ho} and \cite{ct:qs} imply that all elements $\phi^{}_{m\tss a}$
and $\psi^{}_{m\tss a}$, as well as the coefficients of the polynomials
$\tr\ts (\tau+E[-1])^m$
belong to the Feigin--Frenkel center $\z(\wh\gl_N)$; see also
\cite{mr:mm} for a simpler proof. Moreover, each of the families
\ben
\phi^{}_{1\tss 1},\dots,\phi^{}_{N N}\Fand
\psi^{}_{1\tss 1},\dots,\psi^{}_{N N}
\een
is a {\it complete set of Segal--Sugawara vectors\/} in the sense
that the elements of each family together with their images under
all positive powers of the translation
operator $T={\rm ad}\ts\tau$ are algebraically independent and generate
$\z(\wh\gl_N)$.

The elements $\mu_i=E_{ii}$
with $i=1,\dots,N$ span a
Cartan subalgebra of $\gl_N$.
Elements of the classical $\Wc$-algebra $\Wc(\gl_N)$
are regarded as polynomials in the $\mu_i[r]$ with $r<0$.
A calculation of the images of the polynomials \eqref{deftra} with $m=N$
and $\tr\ts (\tau+E[-1])^m$ under the isomorphism \eqref{hchiaffe}
was given in \cite{cm:ho}. The same method applies to
all polynomials \eqref{deftra} and \eqref{deftrh} to
yield the formulas
\begin{align}\label{hcha}
\chi:\tr\ts A^{(m)}\big(\tau+E[-1]_1\big)\dots&\big(\tau+E[-1]_m\big)
\mapsto e_m\big(\tau+\mu^{}_1[-1],\dots,\tau+\mu^{}_N[-1]\big),\\[0.7em]
\label{hchh}
\chi:\tr\ts H^{(m)}\big(\tau+E[-1]_1\big)\dots&\big(\tau+E[-1]_m\big)
\mapsto h_m\big(\tau+\mu^{}_1[-1],\dots,\tau+\mu^{}_N[-1]\big),
\end{align}
where we use standard noncommutative versions of the complete
and elementary symmetric functions in the ordered variables $x_1,\dots,x_p$
defined by the respective formulas
\begin{align}\label{comp}
h_m(x_1,\dots,x_p)&=\sum_{i_1\leqslant\dots\leqslant i_m}
x_{i_1}\dots x_{i_m},\\
\label{elem}
e_m(x_1,\dots,x_p)&=\sum_{i_1>\dots> i_m}
x_{i_1}\dots x_{i_m}.
\end{align}
Relations \eqref{hcha} and \eqref{hchh} can also be derived
from the Yangian character formulas as we indicate below; see Secs~\ref{subsec:csgln}
and \ref{sec:gw}.

\medskip
\noindent
{\bf 1.3.}\quad
Now turn to the Lie algebras of types $B$, $C$ and $D$ and
let $\g=\g_N$ be the orthogonal Lie algebra $\oa_N$
(with $N=2n$ or $N=2n+1$) or
the symplectic Lie algebra $\spa_N$ (with $N=2n$).
We will use the elements
$F_{ij}[r]$ of the loop algebra $\g_N[t,t^{-1}]$, where the $F_{ij}$
are standard generators of $\g_N$;
see Sec.~\ref{subsec:yagn} for the definitions.
For each $a\in\{1,\dots,m\}$
introduce the element $F[r]_a$ of the algebra \eqref{tenprka}
by
\beql{matnot}
F[r]_a=\sum_{i,j=1}^{N}
1^{\ot(a-1)}\ot e_{ij}\ot 1^{\ot(m-a)}\ot F_{ij}[r],
\eeq
where $\U$ in \eqref{tenprka} now stands
for the universal enveloping algebra of
$\wh\g_N\oplus\CC\tau$.
We let $S^{(m)}$ denote the
element of the algebra \eqref{tenprka} which is the image of the
symmetrizer of the Brauer algebra $\Bc_m(\om)$ under
its natural action on $(\CC^{N})^{\ot m}$, where the parameter $\om$
should be specialized to $N$ or $-N$ in the orthogonal and symplectic case, respectively.
The component
of $S^{(m)}$ in $\U$ is the identity; see \eqref{symo} and \eqref{symsp}
below for explicit formulas.
We will use the notation
\beql{galm}
\ga_m(\om)=\frac{\om+m-2}{\om+2\tss m-2}
\eeq
and define the elements $\phi^{}_{m\tss a}\in\U\big(t^{-1}\g_N[t^{-1}]\big)$
by the expansion
\beql{deftr}
\ga_m(\om)\ts\tr\ts S^{(m)} \big(\tau+F[-1]_1\big)\dots \big(\tau+F[-1]_m\big)
=\phi^{}_{m\tss0}\ts\tau^m+\phi^{}_{m\tss1}\ts\tau^{m-1}
+\dots+\phi^{}_{m\tss m},
\eeq
where the trace is taken
over all $m$ copies of $\End\CC^{N}$ (we included the constant factor \eqref{galm}
to get a uniform expression in all cases).
By the main result of \cite{m:ff}, all coefficients $\phi^{}_{m\tss a}$
belong to the Feigin--Frenkel center $\z(\wh\g_N)$.
Note that in the symplectic case $\g_N=\spa_{2n}$
the values of $m$ were restricted to $1\leqslant m\leqslant 2n$,
but the result and arguments also extend to $m=2n+1$; see \cite[Sec.~3.3]{m:ff}.
In the even orthogonal case $\g_N=\oa_{2n}$ there is an additional element
$\phi^{\tss\prime}_n=\Pf\ts\wt F[-1]$ of the center defined as
the (noncommutative) Pfaffian of the skew-symmetric matrix $\wt F[-1]=[\wt F_{ij}[-1]]$,
\beql{genpf}
\Pf\ts\wt F[-1]=\frac{1}{2^nn!}\sum_{\si\in\Sym_{2n}}\sgn\si\cdot
\wt F_{\si(1)\ts\si(2)}[-1]\dots \wt F_{\si(2n-1)\ts\si(2n)}[-1],
\eeq
where $\wt F_{ij}[-1]=F_{ij'}[-1]$ with
$i'=2\tss n-i+1$. The family
$\phi^{}_{2\tss 2},\phi^{}_{4\tss 4},\dots,
\phi^{}_{2n\ts 2n}$
is a complete set of Segal--Sugawara vectors
for $\oa_{2n+1}$ and $\spa_{2n}$, while
$\phi^{}_{2\tss 2},\phi^{}_{4\tss 4},\dots,\phi^{}_{2n-2\ts 2n-2},
\phi^{\tss\prime}_n$
is a complete set of Segal--Sugawara vectors
for $\oa_{2n}$.

The Lie algebras $\oa_{2n+1}$ and $\spa_{2n}$ are Langlands dual to each other,
while $\oa_{2n}$ is self-dual. In all the cases we denote by $\h$
the Cartan subalgebra of $\g_N$ spanned by the elements $\mu_i=F_{ii}$ with
$i=1,\dots,n$
and identify it with
the Cartan subalgebra of ${}^L\g_N$ spanned by the elements with the same names.
We let $\mu_i[r]=\mu_i\ts t^r$ with $r<0$
and $i=1,\dots,n$ denote the basis elements of the vector space
$t^{-1}\h[t^{-1}]$ so that the elements of the classical $\Wc$-algebra $\Wc({}^L\g_N)$
are regarded as polynomials in the $\mu_i[r]$.

\begin{mthm}
The image of the polynomial \eqref{deftr} under
the isomorphism \eqref{hchiaffe}
is given by the formula{\tss\rm:}
\ben
\bal
&\text{type $B_n${\rm:}\qquad\qquad} h_m\big(\tau+\mu_{1}[-1],\dots,
\tau+\mu_{n}[-1],\tau-\mu_{n}[-1],\dots
\tau-\mu_{1}[-1]\big),\\[1.5em]
&\text{type $D_n${\rm:}\qquad\qquad}{\textstyle \frac{1}{2}}\ts h_m\big(\tau+\mu_{1}[-1],\dots,
\tau+\mu_{n-1}[-1],\tau-\mu_{n}[-1],\dots
\tau-\mu_{1}[-1]\big)\\[0.5em]
{}&{\qquad\qquad\qquad\qquad}\quad+{\textstyle \frac{1}{2}}\ts h_m\big(\tau+\mu_{1}[-1],\dots,
\tau+\mu_{n}[-1],\tau-\mu_{n-1}[-1],\dots
\tau-\mu_{1}[-1]\big),\\[1.5em]
&\text{type $C_n${\rm:}\qquad\qquad}e_m\big(\tau+\mu_{1}[-1],\dots,
\tau+\mu_{n}[-1],\tau,\tau-\mu_{n}[-1],\dots
\tau-\mu_{1}[-1]\big).
\eal
\een
Moreover, the image
of the element $\phi^{\tss\prime}_n$ in type $D_n$ is given by
\beql{pfaim}
\big(\mu_{1}[-1]-\tau\big)\dots \big(\mu_{n}[-1]-\tau\big)\ts 1.
\eeq
\end{mthm}

In the last relation $\tau$ is understood as the differentiation operator
so that $\tau\ts 1=0$.

\medskip
\noindent
{\bf 1.4.}\quad
The Fourier coefficients of the image of
any element of the Feigin--Frenkel center $\z(\wh\g)$ under
the state-field correspondence map are well-defined operators
(called the {\it Sugawara operators\/})
on the Wakimoto modules over $\wh\g$. These operators act by
multiplication by scalars which are determined by
the Harish-Chandra image under the isomorphism \eqref{hchiaff};
see \cite[Ch.~8]{f:lc}.
Therefore, the Main Theorem yields explicit formulas for the eigenvalues
of a family of the (higher) Sugawara operators
in the Wakimoto modules.

Our approach is based on the theory of characters originated
in \cite{k:st} in the Yangian context and
in \cite{fr:qc} in the context of quantum affine algebras
(the latter are commonly known as the $q$-characters). The theory
was further developed in \cite{fm:cq} where an algorithm
for the calculation of the $q$-characters was proposed, while
conjectures for functional relations satisfied by the $q$-characters
were proved in \cite{h:kr} and \cite{n:ta}.
An extensive review of the role of the
$q$-characters in classical and quantum integrable systems
is given in \cite{kns:ts}; see also earlier papers \cite{kos:qj}, \cite{ks:ab},
\cite{nn:pt} and \cite{nn:pd} where some
formulas concerning the representations we dealing with in this paper
had been conjectured and studied. In recent work \cite{my:pd},
\cite{my:et} the $q$-characters have been calculated for a wide class
of representations in type $B$, and associated {\it extended $T$-systems\/}
have been introduced.

Due to the general results on the connection of the $q$-characters
with the Feigin--Frenkel center and the classical $\Wc$-algebras
described in \cite[Sec.~8.5]{fr:qc}, one could expect that the
character formulas would be useful for the calculation of the
Harish-Chandra images of the coefficients of the polynomial
\eqref{deftr}. Indeed, as we demonstrate below,
the images in the classical $\Wc$-algebra are closely related with the top
degree components of some linear combinations of the $q$-characters.

We now briefly describe the contents of the paper. We start by
proving the character formulas for some classes of representations
of the Yangian $\Y(\g_N)$ associated with the Lie algebra $\g_N$
(Sec.~\ref{sec:cy}). To this end we employ realizations
of the representations in harmonic tensors
and construct special bases of the representation spaces.
The main calculation is given in Sec.~\ref{sec:ic}, where we consider particular
linear combinations of the Yangian characters and calculate their top degree terms
as elements of the associated graded algebra $\gr\ts\Y(\g_N)\cong \U(\g_N[t])$.
In Sec.~\ref{sec:cw} we recall the definition of the classical
$\Wc$-algebras and write explicit screening operators in all classical types.
By translating the results of Sec.~\ref{sec:ic} to the universal enveloping algebra
$\U(t^{-1}\g_N[t^{-1}])$ we will be able to get them in the form
provided by the Main Theorem (Sec.~\ref{sec:gw}). Finally, in Sec~\ref{sec:ce}
we apply our results to get the Harish-Chandra images
of the Casimir elements for
the Lie algebras $\g_N$ arising from the Brauer--Schur--Weyl duality.
We show that our formulas are equivalent to those
previously found in \cite{imr:ce}.

\section{Characters of Yangian representations}
\label{sec:cy}
\setcounter{equation}{0}

\subsection{Yangian for $\gl_N$}
\label{subsec:yagln}

Denote by $\h$ the Cartan subalgebra of $\gl_N$ spanned by the basis
elements $E_{11},\dots,E_{NN}$. The highest weights of representations
of $\gl_N$ will be considered with respect to this basis,
and the highest vectors will be assumed to be annihilated by the action of the elements
$E_{ij}$ with $1\leqslant i<j\leqslant N$, unless stated otherwise.

Recall the $RTT$-presentation of the Yangian
associated with the Lie algebra $\gl_N$; see e.g. \cite[Ch.~1]{m:yc}.
For $1\leqslant a<b\leqslant m$
introduce the elements $P_{a\tss b}$ of the tensor product algebra
\beql{tenprke}
\underbrace{\End\CC^{N}\ot\dots\ot\End\CC^{N}}_m
\eeq
by
\beql{pdef}
P_{a\tss b}=\sum_{i,j=1}^N 1^{\ot(a-1)}\ot e_{ij}
\ot 1^{\ot(b-a-1)}\ot e_{ji}\ot 1^{\ot(m-b)}.
\eeq
The Yang $R$-matrix $R_{12}(u)$ is a rational function in a complex parameter $u$
with values in the tensor product algebra
$\End\CC^N\ot\End\CC^N$ defined by
\ben
R_{12}(u)=1-\frac{P_{12}}{u}.
\een
This function satisfies the Yang--Baxter equation
\beql{yberep}
R_{12}(u)\ts R_{13}(u+v)\ts R_{23}(v)
=R_{23}(v)\ts R_{13}(u+v)\ts R_{12}(u),
\eeq
where the subscripts indicate the copies of $\End\CC^N$ in the algebra
\eqref{tenprke} with $m=3$.
The {\it Yangian\/}
$\Y(\gl_N)$
is an associative algebra with generators
$t_{ij}^{(r)}$, where $1\leqslant i,j\leqslant N$ and $r=1,2,\dots$,
satisfying certain quadratic relations. To write them down,
introduce the formal series
\ben
t_{ij}(u)=\de_{ij}+\sum_{r=1}^{\infty}t_{ij}^{(r)}\ts u^{-r}
\in\Y(\gl_N)[[u^{-1}]]
\een
and set
\ben
T(u)=\sum_{i,j=1}^N e_{ij}\ot t_{ij}(u)
\in \End\CC^N\ot \Y(\gl_N)[[u^{-1}]].
\een
Consider the algebra
$\End\CC^N\ot\End\CC^N\ot \Y(\gl_N)[[u^{-1}]]$
and introduce its elements $T_1(u)$ and $T_2(u)$ by
\beql{T1T2}
T_1(u)=\sum_{i,j=1}^N e_{ij}\ot 1\ot t_{ij}(u),\qquad
T_2(u)=\sum_{i,j=1}^N 1\ot e_{ij}\ot t_{ij}(u).
\eeq
The defining relations for the algebra $\Y(\gl_N)$ can then
be written in the form
\beql{RTT}
R_{12}(u-v)\ts T_1(u)\ts T_2(v)=T_2(v)\ts T_1(u)\ts R_{12}(u-v).
\eeq

We identify the universal enveloping algebra $\U(\gl_N)$
with a subalgebra of the Yangian $\Y(\gl_N)$ via the embedding
$E_{ij}\mapsto t_{ij}^{(1)}$. Then $\Y(\gl_N)$ can be regarded as a $\gl_N$-module
with the adjoint action. Denote by $\Y(\gl_N)^{\h}$ the
subalgebra of $\h$-invariants under this action.
Consider the left ideal $I$ of the algebra $\Y(\gl_N)$ generated by
all elements $t^{(r)}_{ij}$ with the conditions
$1\leqslant i<j\leqslant N$ and $r\geqslant 1$.
By the Poincar\'e--Birkhoff--Witt
theorem for the Yangian \cite[Sec.~1.4]{m:yc}, the intersection
$\Y(\gl_N)^{\h}\cap I$ is a two-sided ideal of $\Y(\gl_N)^{\h}$. Moreover,
the quotient of $\Y(\gl_N)^{\h}$ by this ideal is isomorphic
to the commutative algebra freely generated by the images
of the elements $t_{ii}^{(r)}$ with $i=1,\dots,N$ and $r\geqslant 1$
in the quotient. We will use the notation $\la^{(r)}_i$
for this image of $t_{ii}^{(r)}$. Thus, we get an analogue of the Harish-Chandra
homomorphism \eqref{hchcl},
\beql{yhch}
\Y(\gl_N)^{\h}\to\CC[\la^{(r)}_i\ts|\ts i=1,\dots,N,\ r\geqslant 1].
\eeq
We combine the elements $\la^{(r)}_i$ into the formal series
\ben
\la_i(u)=1+\sum_{r=1}^{\infty} \la^{(r)}_i\ts u^{-r},\qquad i=1,\dots,N,
\een
which can be understood as the images of the series
$t_{ii}(u)$ under the homomorphism \eqref{yhch}.

The {\it symmetrizer\/} $H^{(m)}$ and {\it anti-symmetrizer\/} $A^{(m)}$
in the algebra \eqref{tenprke}
are the operators in the tensor product space $(\CC^N)^{\ot m}$
associated with the corresponding idempotents
in the group algebra of the symmetric group $\Sym_m$
via its natural
action on
the tensor product space $(\CC^N)^{\ot m}$.
That is,
\beql{symasym}
H^{(m)}=\frac{1}{m!}\ts\sum_{s\in \Sym_m} P_s
\Fand
A^{(m)}=\frac{1}{m!}\ts\sum_{s\in \Sym_m} \sgn s\cdot P_s,
\eeq
where $P_s$ is the element of the algebra \eqref{tenprke}
corresponding to $s\in \Sym_m$.
Both the symmetrizer and anti-symmetrizer admit
multiplicative expressions in terms of the values of the Yang $R$-matrix,
\beql{raprah}
H^{(m)}=\frac{1}{m!}
\prod_{1\leqslant a<b\leqslant m}
\Big(1+\frac{P_{a\tss b}}{b-a}\Big)\Fand
A^{(m)}=\frac{1}{m!}
\prod_{1\leqslant a<b\leqslant m}
\Big(1-\frac{P_{a\tss b}}{b-a}\Big),
\eeq
where the products are taken in the lexicographic order
on the pairs $(a,b)$; see e.g. \cite[Sec.~6.4]{m:yc}.
The operators $H^{(m)}$ and $A^{(m)}$ project $(\CC^N)^{\ot m}$ to the subspaces
of symmetric and skew-symmetric tensors, respectively. Both subspaces
carry irreducible representations of the Yangian $\Y(\gl_N)$.
Consider the tensor product algebra
\beql{tenprkya}
\underbrace{\End\CC^{N}\ot\dots\ot\End\CC^{N}}_m{}\ot \Y(\gl_N)[[u^{-1}]]
\eeq
and extend the notation \eqref{T1T2} to elements of \eqref{tenprkya}.
All coefficients of the formal series
\begin{align}\label{betheh}
&\tr\ts H^{(m)}T_1(u)\ts T_2(u+1)\dots T_m(u+m-1)\\[-1em]
\intertext{and}
&\tr\ts A^{(m)}T_1(u)\ts T_2(u-1)\dots T_m(u-m+1)
\label{bethea}
\end{align}
belong to a commutative subalgebra of the Yangian. This subalgebra
is contained in $\Y(\gl_N)^{\h}$.
The next proposition is well-known and easy to prove;
see also \cite[Sec.~7.4]{bk:rs}, \cite[Sec.~4.5]{fm:ha} and
\cite[Sec.~8.5]{m:yc} for derivations of more general formulas
for the characters of the evaluation modules over $\Y(\gl_N)$.
We give a proof of the proposition to
stress the similarity of the approaches for
all classical types.

\bpr\label{prop:hcha}
The images of the series \eqref{betheh} and \eqref{bethea}
under the homomorphism \eqref{yhch}
are given by
\begin{align}\label{imbethea}
&\sum_{1\leqslant i_1\leqslant\dots\leqslant i_m\leqslant N}
\la_{i_1}(u)\ts\la_{i_2}(u+1)\dots \la_{i_m}(u+m-1)\\[-1em]
\intertext{and}
&\sum_{1\leqslant i_1<\dots< i_m\leqslant N}
\la_{i_1}(u)\ts\la_{i_2}(u-1)\dots \la_{i_m}(u-m+1),
\label{imbetheh}
\end{align}
respectively.
\epr

\bpf
By relations \eqref{RTT} and \eqref{raprah} we can write the product
occurring in \eqref{betheh} as
\beql{stth}
H^{(m)}T_1(u)\dots T_m(u+m-1)=
T_m(u+m-1)\dots T_1(u)\ts H^{(m)}.
\eeq
This relation
shows that the product on
each side can be regarded as
an operator on $(\CC^N)^{\ot m}$ with coefficients in the algebra $\Y(\gl_N)[[u^{-1}]]$
such that the subspace $H^{(m)}\tss (\CC^N)^{\ot m}$ is invariant under this operator.
A basis of this subspace is comprised by vectors of the form
$v_{\ts i_1,\dots,i_m}=H^{(m)}\tss(e_{i_1}\ot\dots\ot e_{i_m})$,
where $i_1\leqslant\dots\leqslant i_m$ and
$e_1,\dots,e_N$ denote the canonical basis vectors of $\CC^N$.
To calculate the trace of the operator, we will find the diagonal
matrix elements corresponding to the basis vectors.
Applying
the operator which occurs on the right hand side of \eqref{stth}
to a basis vector $v_{\ts i_1,\dots,i_m}$ we get
\ben
T_m(u+m-1)\dots T_1(u)\tss H^{(m)}\tss v_{\ts i_1,\dots,i_m}
=T_m(u+m-1)\dots T_1(u)\tss  v_{\ts i_1,\dots,i_m}.
\een
The coefficient of $v_{\ts i_1,\dots,i_m}$ in the expansion of
this expression as a linear combination of the basis vectors is determined
by the coefficient of the tensor $e_{i_1}\ot\dots\ot e_{i_m}$.
Hence, a nonzero contribution to the
image of the diagonal matrix element corresponding
to $v_{\ts i_1,\dots,i_m}$ under the homomorphism \eqref{yhch} only comes from
the term  $t_{i_mi_m}(u+m-1)\dots t_{i_1i_1}(u)$.
The sum over all basis vectors yields the resulting
formula for the image of the element \eqref{betheh}.

The calculation of the image of the series \eqref{bethea} is quite similar.
It relies on the identity
\ben
A^{(m)}T_1(u)\dots T_m(u-m+1)=
T_m(u-m+1)\dots T_1(u)\ts A^{(m)}
\een
and a calculation of the diagonal matrix elements of the operator
which occurs on the right hand side on
the basis vectors $A^{(m)}\tss(e_{i_1}\ot\dots\ot e_{i_m})$,
where $i_1<\dots< i_m$.
\epf

\subsection{Yangians for $\oa_N$ and $\spa_N$}
\label{subsec:yagn}

Throughout the paper we use the involution
on the set $\{1,\dots,N\}$ defined by $i'=N-i+1$.
The Lie subalgebra
of $\gl_N$ spanned by the elements $F_{ij}=E_{ij}-E_{j'i'}$
with $i,j\in\{1,\dots,N\}$
is isomorphic to
the orthogonal Lie algebra $\oa_N$. Similarly, the Lie subalgebra
of $\gl_{2n}$ spanned by the elements
$F_{ij}=E_{ij}-\ve_i\tss\ve_j\tss E_{j'i'}$ with
$i,j\in\{1,\dots,2n\}$
is isomorphic to
the symplectic Lie algebra $\spa_{2n}$, where $\ve_i=1$ for $i=1,\dots,n$ and
$\ve_i=-1$ for $i=n+1,\dots,2n$.
We will keep the notation $\g_N$ for the Lie algebra $\oa_N$
(with $N=2n$ or $N=2n+1$) or $\spa_N$ (with $N=2n$).
Denote by $\h$ the Cartan subalgebra of $\g_N$ spanned by the basis
elements $F_{11},\dots,F_{nn}$. The highest weights of representations
of $\g_N$ will be considered with respect to this basis,
and the highest vectors will be assumed to be annihilated by the action of the elements
$F_{ij}$ with $1\leqslant i<j\leqslant N$, unless stated otherwise.

Recall the $RTT$-presentation of the Yangian
associated with the Lie algebra $\g_N$ following
the general approach of \cite{d:qg} and \cite{rtf:ql}; see also \cite{aacfr:rp}
and \cite{amr:rp}.

For $1\leqslant a<b\leqslant m$ consider
the elements $P_{a\tss b}$ of the tensor product algebra
\eqref{tenprke} defined by \eqref{pdef}. Introduce also
the elements $Q_{a\tss b}$ of \eqref{tenprke}
which are defined by different formulas in the orthogonal
and symplectic cases. In the orthogonal case we set
\ben
Q_{a\tss b}=\sum_{i,j=1}^N 1^{\ot(a-1)}\ot e_{ij}
\ot 1^{\ot(b-a-1)}\ot e_{i'j'}\ot 1^{\ot(m-b)},
\een
and in the symplectic case
\ben
Q_{a\tss b}=\sum_{i,j=1}^N \ve_i\tss\ve_j\ts 1^{\ot(a-1)}\ot e_{ij}
\ot 1^{\ot(b-a-1)}\ot e_{i'j'}\ot 1^{\ot(m-b)}.
\een

Set $\kappa=N/2-1$ in the orthogonal case and
$\kappa=N/2+1$ in the symplectic case.
The $R$-matrix $R_{12}(u)$ is a rational function in a complex parameter $u$
with values in the tensor product algebra
$\End\CC^N\ot\End\CC^N$ defined by
\ben
R_{12}(u)=1-\frac{P_{12}}{u}+\frac{Q_{12}}{u-\kappa}.
\een
It is well known by \cite{zz:rf}
that this function satisfies the Yang--Baxter equation \eqref{yberep}.

The {\it Yangian\/}
$\Y(\g_N)$
is an associative algebra with generators
$t_{ij}^{(r)}$, where $1\leqslant i,j\leqslant N$ and $r=1,2,\dots$,
satisfying certain quadratic relations.
Introduce the formal series
\ben
t_{ij}(u)=\de_{ij}+\sum_{r=1}^{\infty}t_{ij}^{(r)}\ts u^{-r}
\in\Y(\g_N)[[u^{-1}]]
\een
and set
\ben
T(u)=\sum_{i,j=1}^N e_{ij}\ot t_{ij}(u)
\in \End\CC^N\ot \Y(\g_N)[[u^{-1}]].
\een
Consider the algebra
$\End\CC^N\ot\End\CC^N\ot \Y(\g_N)[[u^{-1}]]$
and introduce its elements $T_1(u)$ and $T_2(u)$ by
the same formulas \eqref{T1T2} as in the case of $\gl_N$.
The defining relations for the algebra $\Y(\g_N)$ can then
be written in the form
\beql{RTTbcd}
R_{12}(u-v)\ts T_1(u)\ts T_2(v)=T_2(v)\ts T_1(u)\ts R_{12}(u-v)
\eeq
together with the relation
\ben
T^{\tss\prime}(u+\ka)\ts T(u)=1,
\een
where the prime denotes the matrix transposition defined
for an $N\times N$ matrix $A=[A_{ij}]$ by
\ben
(A')_{ij}=A_{j'i'}\Fand (A')_{ij}=\ve_i\tss\ve_j\ts A_{j'i'}
\een
in the orthogonal and symplectic case, respectively.

We identify the universal enveloping algebra $\U(\g_N)$
with a subalgebra of the Yangian $\Y(\g_N)$ via the embedding
\ben
F_{ij}\mapsto t_{ij}^{(1)},\qquad i,j=1,\dots,N.
\een
Then $\Y(\g_N)$ can be regarded as a $\g_N$-module
with the adjoint action. Denote by $\Y(\g_N)^{\h}$ the
subalgebra of $\h$-invariants under this action.

Consider the left ideal $I$ of the algebra $\Y(\g_N)$ generated by
all elements $t^{(r)}_{ij}$ with the conditions
$1\leqslant i<j\leqslant N$ and $r\geqslant 1$.
It follows from the Poincar\'e--Birkhoff--Witt
theorem for the Yangian \cite[Sec.~3]{amr:rp} that the intersection
$\Y(\g_N)^{\h}\cap I$ is a two-sided ideal of $\Y(\g_N)^{\h}$. Moreover,
the quotient of $\Y(\g_N)^{\h}$ by this ideal is isomorphic
to the commutative algebra freely generated by the images
of the elements $t_{ii}^{(r)}$ with $i=1,\dots,n$ and $r\geqslant 1$
in the quotient. We will use the notation $\la^{(r)}_i$
for this image of $t_{ii}^{(r)}$ and extend this notation to
all values $i=1,\dots,N$. Thus, we get an analogue of the Harish-Chandra
homomorphism \eqref{hchcl},
\beql{yhchbcd}
\Y(\g_N)^{\h}\to\CC[\la^{(r)}_i\ts|\ts i=1,\dots,n,\ r\geqslant 1].
\eeq
We combine the elements $\la^{(r)}_i$ into the formal series
\ben
\la_i(u)=1+\sum_{r=1}^{\infty} \la^{(r)}_i\ts u^{-r},\qquad i=1,\dots,N
\een
which can be understood as the image of the series
$t_{ii}(u)$ under the homomorphism \eqref{yhchbcd}.

It follows from \cite[Prop.~5.2 and 5.14]{amr:rp}, that the series $\la_i(u)$
satisfy the relations
\beql{larel}
\la_i(u+\ka-i)\ts\la_{i'}(u)=\la_{i+1}(u+\ka-i)\ts\la_{(i+1)'}(u),
\eeq
for $i=0,1,\dots,n-1$ if $\g_N=\oa_{2n}$ or $\spa_{2n}$,
and for $i=0,1,\dots,n$ if $\g_N=\oa_{2n+1}$, where $\la_0(u)=\la_{0'}(u):=1$.
Under an appropriate identification, the relations \eqref{larel} coincide
with those for the $q$-characters, as the $\la_i(u)$ correspond to
the ``single box variables"; see for instance
\cite[Sec.~7]{kns:ts} and \cite[Sec.~2]{nn:pt}.
This coincidence is consistent with the general result which establishes
the equivalence of the definitions of $q$-characters in \cite{fr:qc}
and \cite{k:st}; see \cite[Prop.~2.4]{fm:cq} for a proof.
The $q$-characters
have been extensively studied; see \cite{fm:ha},
\cite{fr:qc} and \cite{k:st}. In
particular, formulas for the $q$-characters
of some classes of modules were conjectured in
\cite{kos:qj}, \cite{nn:pt}
and \cite{nn:pd} and later proved in \cite{h:kr} and \cite{n:ta}.
However, this was done in the
context of the new realization of the quantum affine algebras. In what follows we
compute some $q$-characters independently in our setting of the $RTT$
realization of the Yangians.

Introduce
the element $S^{(m)}$ of the algebra \eqref{tenprke} by setting $S^{(1)}=1$
and for $m\geqslant 2$ define it
by the respective formulas in the orthogonal and symplectic cases:
\begin{align}\label{symo}
S^{(m)}&=\frac{1}{m!}
\prod_{1\leqslant a<b\leqslant m}
\Big(1+\frac{P_{a\tss b}}{b-a}-\frac{Q_{a\tss b}}
{N/2+b-a-1}\Big)\\[-1.5em]
\intertext{and}
\label{symsp}
S^{(m)}&=\frac{1}{m!}
\prod_{1\leqslant a<b\leqslant m}
\Big(1-\frac{P_{a\tss b}}{b-a}-\frac{Q_{a\tss b}}
{n-b+a+1}\Big),
\end{align}
where the products are taken in the lexicographic order
on the pairs $(a,b)$ and the condition $m\leqslant n+1$
is assumed in \eqref{symsp}.
The elements \eqref{symo} and \eqref{symsp} are the images of the
symmetrizers in the corresponding Brauer algebras $\Bc_m(N)$ and $\Bc_m(-N)$
under their actions on the vector space $(\CC^N)^{\ot m}$. In particular,
for any $1\leqslant a<b\leqslant m$ for the operator $S^{(m)}$ we have
\beql{sqpo}
S^{(m)}\ts Q_{a\tss b}=Q_{a\tss b}\ts S^{(m)}=0
\Fand S^{(m)}\ts P_{a\tss b}=P_{a\tss b}\ts S^{(m)}=\pm S^{(m)}
\eeq
with the plus and minus signs taken in the orthogonal and symplectic
case, respectively. The symmetrizer admits a few
other equivalent expressions which are reproduced in \cite{m:ff}.

In the orthogonal case
the operator $S^{(m)}$ projects $(\CC^N)^{\ot m}$ to the irreducible
representation of the Lie algebra $\oa_N$ with the highest weight $(m,0,\dots,0)$.
The dimension of this representation equals
\ben
\frac{N+2\tss m-2}{N+m-2}\ts\binom{N+m-2}{m}.
\een
This representation is extended to the Yangian $\Y(\oa_N)$
and it is one of the Kirillov--Reshetikhin modules.
In the symplectic case with $m\leqslant n$ the operator
$S^{(m)}$ projects $(\CC^{2n})^{\ot m}$ to
the subspace of skew-symmetric harmonic tensors which carries
an irreducible representation of $\spa_{2n}$ with
the highest weight $(1,\dots,1,0,\dots,0)$
(with $m$ copies of $1$). Its dimension equals
\beql{dimhm}
\frac{2\tss n-2\tss m+2}{2\tss n-m+2}\ts\binom{2\tss n+1}{m}=\binom{2\tss n}{m}
-\binom{2\tss n}{m-2}.
\eeq
This representation is extended to the $m$-th
fundamental representation of the Yangian $\Y(\spa_{2n})$
which is also a Kirillov--Reshetikhin module.
It is well-known that if $m=n+1$ then
the subspace of tensors is zero
so that $S^{(n+1)}=0$.

The
existence of the Yangian action on the Lie algebra modules here
can be explained by the fact that the projections \eqref{symo} and \eqref{symsp}
are the products of the evaluated $R$-matrices
\beql{rmapr}
S^{(m)}=\frac{1}{m!}
\prod_{1\leqslant a<b\leqslant m} R_{a\tss b}(u_a-u_b),
\eeq
where $u_a=u+a-1$ and $u_a=u-a+1$ for $a=1,\dots,m$ in the orthogonal and
symplectic case, respectively; see \cite{imo:nf} for a proof
in the context of a fusion procedure for the Brauer algebra.
The same fact leads
to a construction of a commutative subalgebra of the Yangian $\Y(\g_N)$;
see \cite{m:ff}. We will calculate the images of the elements of this
subalgebra under the homomorphism \eqref{yhchbcd} and thus reproduce
the character formulas for the respective classes of Yangian representations;
cf. \cite[Sec.~7]{kns:ts}.
Consider the tensor product algebra
\beql{tenprky}
\underbrace{\End\CC^{N}\ot\dots\ot\End\CC^{N}}_m{}\ot \Y(\g_N)[[u^{-1}]]
\eeq
and extend the notation \eqref{T1T2} to elements of \eqref{tenprky}.

\subsubsection{Series $B_n$}

The commutative subalgebra of the Yangian $\Y(\oa_{N})$ with $N=2n+1$ is generated by the
coefficients of the formal series
\beql{betheb}
\tr\ts S^{(m)}T_1(u)\ts T_2(u+1)\dots T_m(u+m-1)
\eeq
with the trace taken over all $m$ copies of $\End\CC^{N}$ in \eqref{tenprky},
where $\g_N=\oa_{N}$ and $S^{(m)}$ is defined in \eqref{symo}. It follows easily
from the defining relations \eqref{RTTbcd} that all elements of this subalgebra
belong to $\Y(\oa_N)^{\h}$.

\bpr\label{prop:hchb}
The image of the series \eqref{betheb} under the homomorphism \eqref{yhchbcd}
is given by
\ben
\sum_{1\leqslant i_1\leqslant\dots\leqslant i_m\leqslant N}
\la_{i_1}(u)\ts\la_{i_2}(u+1)\dots \la_{i_m}(u+m-1)
\een
with the condition that $n+1$ occurs among the summation indices $i_1,\dots,i_m$
at most once.
\epr

\bpf
By \cite[Prop.~3.1]{m:ff} the operator $S^{(m)}$ can be given by the formula
\beql{symmeo}
S^{(m)}=H^{(m)}\ts\sum_{r=0}^{\lfloor m/2\rfloor}
\frac{(-1)^r}{2^{\tss r}\ts r!}\binom{N/2+m-2}{r}^{-1}
\sum_{a_i<b_i}Q_{a_1b_1}Q_{a_2b_2}\dots Q_{a_rb_r}
\eeq
with the second sum taken over the (unordered)
sets of disjoint pairs
$\{(a_1,b_1),\dots,(a_r,b_r)\}$ of indices from $\{1,\dots,m\}$. Here
$H^{(m)}$ is the symmetrization operator defined in \eqref{symasym}.
Note that for each $r$ the second sum
in \eqref{symmeo} commutes with
any element $P_s$ and hence commutes with $H^{(m)}$.

Recall that the subspace of harmonic tensors in $(\CC^N)^{\ot m}$ is spanned
by the tensors $v$ with the property $Q_{a\tss b}\ts v=0$ for all $1\leqslant a<b\leqslant m$.
By \eqref{sqpo} the operator $S^{(m)}$ projects $(\CC^N)^{\ot m}$ to a subspace of symmetric
harmonic tensors which we denote by $\Hc_m$. This subspace carries an
irreducible representation of $\oa_N$
with the highest weight $(m,0,\dots,0)$.
Therefore, the trace in \eqref{betheb}
can be calculated over the subspace $\Hc_m$. We will introduce a special basis of
this subspace. We identify the image of the symmetrizer $H^{(m)}$
with the space of homogeneous polynomials of degree $m$ in variables $z_1,\dots,z_N$
via the isomorphism
\beql{symsub}
H^{(m)}(e_{i_1}\ot\dots\ot e_{i_m})\mapsto z_{i_1}\dots z_{i_m}.
\eeq
The subspace
$\Hc_m$
is then identified with the subspace
of harmonic homogeneous polynomials of degree $m$; they belong to the kernel
of the Laplace operator
\ben
\Delta=\sum_{i=1}^n \di_{z_i}\di_{z_{i'}}+\frac12\ts\di_{z_{n+1}}^2.
\een
The basis vectors of $\Hc_m$ will be parameterized by the $N$-tuples
$(k_1,\dots,k_n,\de,l_n,\dots,l_1)$, where the $k_i$ and $l_i$ are arbitrary
nonnegative integers, $\de\in\{0,1\}$ and the sum of all entries is $m$.
Given such a tuple, the corresponding harmonic polynomial is defined by
\beql{basb}
\sum_{a_1,\dots,a_n}\frac{(-2)^{a_1+\dots+a_n}(a_1+\dots+a_n)!
\ts z_{n+1}^{2a_1+\dots+2a_n+\de}}
{a_1!\dots a_n!\ts
(2a_1+\dots+2a_n+\de)!}\prod_{i=1}^n
\frac{z_i^{k_i-a_i}\ts z_{i'}^{l_i-a_i}}{(k_i-a_i)!\ts(l_i-a_i)!},
\eeq
summed over nonnegative integers $a_i$ satisfying $a_i\leqslant\min\{k_i,l_i\}$.
Each polynomial contains a unique monomial (which we call the {\it leading
monomial\/}) where the variable $z_{n+1}$ occurs with
the power not exceeding $1$.
It is straightforward to see that
these polynomials are all harmonic and linearly independent.
Furthermore, a simple calculation shows that the number of the polynomials
coincides with the dimension of the irreducible representation of $\oa_N$
with the highest weight $(m,0,\dots,0)$ and so they form a basis of the
subspace $\Hc_m$.

By relations \eqref{RTTbcd} and \eqref{rmapr} we can write the product
occurring in \eqref{betheb} as
\beql{stt}
S^{(m)}T_1(u)\dots T_m(u+m-1)=
T_m(u+m-1)\dots T_1(u)\ts S^{(m)}.
\eeq
This relation together with \eqref{sqpo}
shows that the product on
each side can be regarded as
an operator on $(\CC^N)^{\ot m}$ with coefficients in the algebra $\Y(\oa_N)[[u^{-1}]]$
such that the subspace $\Hc_m$ is invariant under this operator.
Now fix a basis vector $v\in \Hc_m$ of the form \eqref{basb}.
Denote the operator on the
right hand side of \eqref{stt} by $A$
and consider the coefficient of $v$ in the expansion of $A\tss v$
as a linear combination of the basis vectors.
Use the isomorphism \eqref{symsub} to write the vector $v$
as a linear combination of the tensors
$e_{j_1}\ot\dots\ot e_{j_m}$.
We have $S^{(m)} v=v$, while the matrix elements of the remaining product
are found from the expansion
\begin{multline}
T_m(u+m-1)\dots T_1(u)(e_{j_1}\ot\dots\ot e_{j_m})\\
{}=\sum_{i_1,\dots,i_m}
t_{i_mj_m}(u+m-1)\dots t_{i_1j_1}(u) (e_{i_1}\ot\dots\ot e_{i_m}).
\non
\end{multline}
The coefficient of $v$ in the expansion of $A\tss v$ is uniquely determined
by the coefficient of the tensor $e_{i_1}\ot\dots\ot e_{i_m}$
with $i_1\leqslant\dots\leqslant i_m$ which corresponds to the leading
monomial of $v$ under the isomorphism \eqref{symsub}. It is clear from
formula \eqref{basb} that if a tensor of the form
$e_{j_1}\ot\dots\ot e_{j_m}$ corresponds to
a non-leading monomial, then the matrix element
$t_{i_mj_m}(u+m-1)\dots t_{i_1j_1}(u)$ vanishes under the homomorphism \eqref{yhchbcd}.
Therefore, a nonzero contribution to the
image of the diagonal matrix element of the operator $A$ corresponding
to $v$ under the homomorphism \eqref{yhchbcd} only comes from
the term  $t_{i_mi_m}(u+m-1)\dots t_{i_1i_1}(u)$.
Taking the sum over all basis vectors \eqref{basb} yields the resulting
formula for the image of the element \eqref{betheb}.
\epf

\subsubsection{Series $D_n$}

The commutative subalgebra of the Yangian $\Y(\oa_{N})$ with $N=2n$ is generated by the
coefficients of the formal series defined by the same formula \eqref{betheb},
where the parameter $N$ now takes an even value $2n$.

\bpr\label{prop:hchd}
The image of the series \eqref{betheb} under the homomorphism \eqref{yhchbcd}
is given by
\ben
\sum_{1\leqslant i_1\leqslant\dots\leqslant i_m\leqslant N}
\la_{i_1}(u)\ts\la_{i_2}(u+1)\dots \la_{i_m}(u+m-1)
\een
with the condition that $n$ and $n'$ do not occur simultaneously
among the summation indices $i_1,\dots,i_m$.
\epr

\bpf
As in the proof of Proposition~\ref{prop:hchb}, we use the formula \eqref{symmeo} for
the symmetrizer $S^{(m)}$ and its properties \eqref{sqpo}. Following the argument
of that proof we identify the image $S^{(m)}(\CC^N)^{\ot m}$
with the space $\Hc_m$ of homogeneous harmonic
polynomials of degree $m$ in variables $z_1,\dots,z_N$ via
the isomorphism \eqref{symsub}. This time the harmonic polynomials are annihilated
by the Laplace operator of the form
\ben
\Delta=\sum_{i=1}^n \di_{z_i}\di_{z_{i'}}.
\een
The basis vectors of $\Hc_m$ will be parameterized by the $N$-tuples
$(k_1,\dots,k_n,l_n,\dots,l_1)$, where the $k_i$ and $l_i$ are arbitrary
nonnegative integers, the sum of all entries is $m$ and at least
one of $k_n$ and $l_n$ is zero.
Given such a tuple, the corresponding harmonic polynomial is now defined by
\begin{multline}\label{basd}
\sum_{a_1,\dots,a_{n-1}}\frac{(-1)^{a_1+\dots+a_{n-1}}(a_1+\dots+a_{n-1})!
\ts z_{n}^{a_1+\dots+a_{n-1}+k_n}\ts z_{n'}^{a_1+\dots+a_{n-1}+l_n}}
{a_1!\dots a_{n-1}!\ts
(a_1+\dots+a_{n-1}+k_n)!\ts (a_1+\dots+a_{n-1}+l_n)!}\\
{}\times\prod_{i=1}^{n-1}
\frac{z_i^{k_i-a_i}\ts z_{i'}^{l_i-a_i}}{(k_i-a_i)!\ts(l_i-a_i)!},
\end{multline}
summed over nonnegative integers $a_1,\dots,a_{n-1}$ satisfying $a_i\leqslant\min\{k_i,l_i\}$.
A unique {\it leading
monomial\/} corresponds to the values $a_1=\dots=a_{n-1}=0$.
The argument is now completed in the same way as for Proposition~\ref{prop:hchb}
by considering the diagonal matrix elements
of the operator on right hand side of \eqref{stt}
corresponding to the basis vectors \eqref{basd}.
These coefficients are determined by those of
the leading monomials and their images under the homomorphism \eqref{yhchbcd}
are straightforward to calculate.
\epf

\subsubsection{Series $C_n$}

The commutative subalgebra of the Yangian $\Y(\spa_{N})$ with $N=2n$ is generated by the
coefficients of the formal series
\beql{bethec}
\tr\ts S^{(m)}T_1(u)\ts T_2(u-1)\dots T_m(u-m+1),
\eeq
with the trace taken over all $m$ copies of $\End\CC^{N}$ in \eqref{tenprky}
with $\g_N=\spa_{N}$ and $S^{(m)}$ defined in \eqref{symsp} with $m\leqslant n$.

\bpr\label{prop:hchc}
The image of the series \eqref{bethec} with $m\leqslant n$
under the homomorphism \eqref{yhchbcd}
is given by
\beql{charc}
\sum_{1\leqslant i_1<\dots< i_m\leqslant 2n}
\la_{i_1}(u)\ts\la_{i_2}(u-1)\dots \la_{i_m}(u-m+1)
\eeq
with the condition that if for any $i$ both $i$ and $i'$ occur
among the summation indices as $i=i_r$ and $i'=i_s$
for some $1\leqslant r<s\leqslant m$, then $s-r\leqslant n-i$.
\epr

\bpf
Using again
\cite[Prop.~3.1]{m:ff} we find that the operator $S^{(m)}$
can be given by the formula
\beql{symmesp}
S^{(m)}=A^{(m)}\ts\sum_{r=0}^{\lfloor m/2\rfloor}
\frac{1}{2^{\tss r}\ts r!}\binom{-n+m-2}{r}^{-1}
\sum_{a_i<b_i}Q_{a_1b_1}Q_{a_2b_2}\dots Q_{a_rb_r}
\eeq
with the second sum taken over the (unordered)
sets of disjoint pairs
$\{(a_1,b_1),\dots,(a_r,b_r)\}$ of indices from $\{1,\dots,m\}$. Here
$A^{(m)}$ is the anti-symmetrization operator
defined in \eqref{symasym}.
For each $r$ the second sum
in \eqref{symmesp} commutes with
any element $P_s$ and hence commutes with $A^{(m)}$.

As with the orthogonal case,
the subspace of harmonic tensors in $(\CC^N)^{\ot m}$ is spanned
by the tensors $v$ with the property $Q_{a\tss b}\ts v=0$ for all $1\leqslant a<b\leqslant m$.
The operator $S^{(m)}$ projects $(\CC^N)^{\ot m}$ to a subspace of skew-symmetric
harmonic tensors which we denote by $\Hc_m$.
Hence, the trace in \eqref{bethec}
can be calculated over the subspace $\Hc_m$. We introduce a special basis of
this subspace by identifying the image of the anti-symmetrizer $A^{(m)}$
with the space of homogeneous polynomials of degree $m$ in
the anti-commuting
variables $\ze_1,\dots,\ze_{2n}$
via the isomorphism
\beql{symsuc}
A^{(m)}(e_{i_1}\ot\dots\ot e_{i_m})\mapsto \ze_{i_1}\wedge\dots\wedge \ze_{i_m}.
\eeq
The subspace $\Hc_m$
is then identified with the subspace
of harmonic homogeneous polynomials of degree $m$; they belong to the kernel
of the Laplace operator
\ben
\Delta=\sum_{i=1}^n \di_{i}\wedge \di_{\tss i^{\tss\prime}},
\een
where $\di_i$ denotes the (left) partial derivative over $\ze_i$.

The basis vectors of $\Hc_m$ will be parameterized by the subsets
$\{i_1,\dots,i_m\}$ of the set $\{1,\dots,2n\}$ satisfying
the condition as stated in the proposition, when the elements $i_1,\dots,i_m$
are written in the increasing order.
We will call such subsets
{\it admissible\/}. The number of admissible subsets can be shown to be
given by the formula \eqref{dimhm}, which coincides with
the dimension of $\Hc_m$.
Consider monomials of the form
\beql{monze}
\ze_{a_1}\wedge \ze_{a'_1}\wedge\dots \wedge\ze_{a_k}\wedge \ze_{a'_k}
\wedge \ze_{b_1}\wedge\dots\wedge \ze_{b_l}
\eeq
with $1\leqslant a_1<\dots<a_k\leqslant n$ and $1\leqslant b_1<\dots<b_l\leqslant 2n$,
associated with subsets $\{a_1,a'_1,\dots,a_k,a'_k,b_1,\dots,b_l\}$ of $\{1,\dots,2n\}$
of cardinality $m=2k+l$, where $b_i\ne b'_{j}$
for all $i$ and $j$. We will suppose that the parameters $b_i$ are fixed
and label the monomial \eqref{monze} by the $k$-tuple $(a_1,\dots,a_k)$.
Furthermore, we order the $k$-tuples and the corresponding
monomials lexicographically.

Now let
the subset $\{a_1,a'_1,\dots,a_k,a'_k,b_1,\dots,b_l\}$ be admissible
and suppose that the parameters $a_1,\dots,a_k$ are fixed too.
We will call the corresponding monomial \eqref{monze}
{\it admissible\/}.
Fix $i\in\{1,\dots,k\}$. Let $s$ be the number of the elements $b_j$
of the subset
satisfying $a_i<b_j<a'_i$. By the admissibility condition applied
to $a_i$ and $a'_i$, we have the inequality
$2(k-i)+s<n-a_i$. Therefore, there exist elements $c_i,\dots,c_k$ satisfying
$a_i< c_i<\dots<c_k\leqslant n$ so that none of $c_j$ or $c'_j$
with $j=i,\dots,k$ belongs to the subset
$\{a_1,a'_1,\dots,a_k,a'_k,b_1,\dots,b_l\}$. Taking the consecutive values
$i=k,k-1,\dots,1$ choose the maximum possible
element $c_i$ at each step. Thus, we get
a family of elements $c_1<\dots<c_k$ uniquely determined by the admissible subset.
In particular, $c_i>a_i$ for all $i$.

Note that our condition on the parameters $b_i$ implies that the monomial
$\ze_{b_1}\wedge\dots\wedge \ze_{b_l}$ is annihilated by the operator $\Delta$.
We denote this monomial by $y$ and set
$x_a=\ze_a\wedge \ze_{a'}$ for $a=1,\dots,n$. The vector
\ben
\sum_{p=0}^k(-1)^p \sum_{1\leqslant d_1<\dots<d_p\leqslant k}
x_{a_1}\wedge \dots\wedge \wh x_{a_{d_1}}
\wedge \dots\wedge \wh x_{a_{d_p}}\wedge \dots\wedge x_{a_k}
\wedge x_{c_{d_1}}\wedge\dots\wedge x_{c_{d_p}}\wedge y,
\een
where the hats indicate the factors to be skipped,
is easily seen to belong to the kernel of the operator $\Delta$
so it is an element of the subspace $\Hc_m$.
Furthermore, these vectors parameterized by all admissible subsets
form a basis of $\Hc_m$. Indeed, the vectors are linearly independent
because
the linear combination defining each vector is uniquely determined
by the admissible monomial $x_{a_1}\wedge \dots\wedge x_{a_k}\wedge y$
which precedes all the other monomials occurring in the linear combination
with respect to the lexicographic order.

Note that apart from the minimal admissible monomial
$x_{a_1}\wedge \dots\wedge x_{a_k}\wedge y$, the linear combination defining
a basis vector may contain some other admissible monomials.
By eliminating such additional admissible monomials with the use of
an obvious induction on the lexicographic order, we can produce
another basis of the space $\Hc_m$ parameterized by all admissible subsets
with the property that each basis vector is given by
a linear combination of monomials of the same form as above,
containing a unique admissible monomial.

By relations \eqref{RTTbcd} and \eqref{rmapr} we can write the product
occurring in \eqref{bethec} as
\beql{sttsp}
S^{(m)}T_1(u)\dots T_m(u-m+1)=
T_m(u-m+1)\dots T_1(u)\ts S^{(m)}
\eeq
and complete the argument exactly as in the proof of Proposition~\ref{prop:hchb}.
Indeed, relations \eqref{sqpo} and \eqref{sttsp}
show that the product on
each side can be regarded as
an operator on $(\CC^N)^{\ot m}$ with coefficients in the algebra $\Y(\spa_N)[[u^{-1}]]$
such that the subspace $\Hc_m$ is invariant under this operator.
Denote the operator on the right hand side of \eqref{sttsp}
by $A$ and let $v$ denote the basis vector of $\Hc_m$
corresponding to an admissible subset $\{i_1,\dots,i_m\}$
with $i_1<\dots< i_m$.
The properties of the basis vectors imply that
a nonzero contribution to the
image of the diagonal matrix element of the operator $A$ corresponding
to $v$ under the homomorphism \eqref{yhchbcd} only comes from
the term  $t_{i_mi_m}(u-m+1)\dots t_{i_1i_1}(u)$.
\epf

We will be using an equivalent formula for the expression
\eqref{charc} given in \cite[Prop.~2.4]{kosy:dl}.
The argument there is combinatorial and relies only on the
identities \eqref{larel}. To state the formula from \cite{kosy:dl}
introduce new parameters $\vk_i(u)$ for $i=1,\dots,2n+2$ by
\ben
\vk_i(u)=\la_i(u),\qquad \vk_{2n-i+3}(u)=\la_{2n-i+1}(u)\qquad\text{for}
\quad i=1,\dots,n,
\een
and $\vk_{n+2}(u)=-\vk_{n+1}(u)$, where $\vk_{n+1}(u)$
is a formal series in $u^{-1}$
with constant term $1$ defined by
\ben
\vk_{n+1}(u)\tss\vk_{n+1}(u-1)=\la_n(u)\tss\la_{n'}(u-1).
\een

\bco\label{cor:eqka}
The image of the series \eqref{bethec} with $m\leqslant n$
under the homomorphism \eqref{yhchbcd}
can be written as
\beql{laka}
\sum_{1\leqslant i_1<\dots< i_m\leqslant 2n+2}
\vk_{i_1}(u)\ts\vk_{i_2}(u-1)\dots \vk_{i_m}(u-m+1).
\eeq
Moreover, the expression \eqref{laka} is zero for $m=n+1$.
\qed
\eco

\section{Harish-Chandra images for the current algebras}
\label{sec:ic}
\setcounter{equation}{0}

We will use the character formulas obtained in Sec.~\ref{sec:cy}
to calculate the Harish-Chandra images of elements of certain commutative
subalgebras of $\U\big(\g[t]\big)$ for the simple Lie algebras
$\g$ of all classical types. The results in the case of $\gl_N$ are well-known,
the commutative subalgebras were constructed explicitly in \cite{t:qg};
see also \cite{cm:ho}, \cite{ct:qs}, \cite{mr:mm}, \cite{mtv:be} and \cite{mtv:ba}.

\subsection{Case of $\gl_N$}
\label{subsec:csgln}

Identify the universal enveloping algebra $\U(\gl_N)$
with a subalgebra of $\U\big(\gl_N[t]\big)$ via the embedding
$E_{ij}\mapsto E_{ij}[0]$. Then $\U\big(\gl_N[t]\big)$ can be regarded as a $\gl_N$-module
with the adjoint action. Denote by $\U\big(\gl_N[t]\big)^{\h}$ the
subalgebra of $\h$-invariants under this action.
Consider the left ideal $I$ of the algebra $\U\big(\gl_N[t]\big)$ generated by
all elements $E_{ij}[r]$ with the conditions
$1\leqslant i<j\leqslant N$ and $r\geqslant 0$.
By the Poincar\'e--Birkhoff--Witt
theorem, the intersection
$\U\big(\gl_N[t]\big)^{\h}\cap I$ is a two-sided ideal of
$\U\big(\gl_N[t]\big)^{\h}$. Moreover,
the quotient of $\U\big(\gl_N[t]\big)^{\h}$ by this ideal is isomorphic
to the commutative algebra freely generated by the images
of the elements $E_{ii}[r]$ with $i=1,\dots,N$ and $r\geqslant 0$
in the quotient. We will denote by $\mu_i[r]$
this image of $E_{ii}[r]$.
We get an analogue of the Harish-Chandra
homomorphism \eqref{hchcl},
\beql{yhchcua}
\U\big(\gl_N[t]\big)^{\h}\to\CC[\mu_i[r]\ts|\ts i=1,\dots,N,\ r\geqslant 0].
\eeq
Combine the elements $E_{ij}[r]$ and $\mu_i[r]$ into the formal series
\ben
E_{ij}(u)=\sum_{r=0}^{\infty} E_{ij}[r]\ts u^{-r-1}\Fand
\mu_i(u)=\sum_{r=0}^{\infty} \mu_i[r]\ts u^{-r-1}.
\een
Then $\mu_i(u)$
is understood as the image of the series
$E_{ii}(u)$ under the homomorphism \eqref{yhchcua}.
Consider tensor product algebras
\ben
\underbrace{\End\CC^{N}\ot\dots\ot\End\CC^{N}}_m{}\ot \U\big(\gl_N[t]\big)[[u^{-1},\di_u]]
\een
and use matrix notation as in \eqref{matnota}.

\bpr\label{prop:hchan}
For the images under the Harish-Chandra homomorphism \eqref{yhchcua} we have
\begin{align}\label{hchan}
\tr\ts A^{(m)} \big(\di_u+E_1(u)\big)\dots \big(\di_u+E_m(u)\big)
&\mapsto
e_m\big(\di_u+\mu_{1}(u),\dots,\di_u+\mu_{N}(u)\big),\\[0.7em]
\tr\ts H^{(m)} \big(\di_u+E_1(u)\big)\dots \big(\di_u+E_m(u)\big)
&\mapsto
h_m\big(\di_u+\mu_{1}(u),\dots,\di_u+\mu_{N}(u)\big).
\label{hchhn}
\end{align}
\epr

\bpf
The argument is essentially the same as in the proof
of Proposition~\ref{prop:hcha}.
Both relations are immediate from the cyclic
property of trace and the identities
\ben
\bal
\big(\di_u+E_1(u)\big)\dots \big(\di_u+E_m(u)\big)A^{(m)}& =
A^{(m)} \big(\di_u+E_1(u)\big)\dots \big(\di_u+E_m(u)\big)A^{(m)},\\[0.7em]
H^{(m)}\big(\di_u+E_1(u)\big)\dots \big(\di_u+E_m(u)\big)& =
H^{(m)} \big(\di_u+E_1(u)\big)\dots \big(\di_u+E_m(u)\big)H^{(m)},
\eal
\een
implied by the fact that $\di_u+E(u)$ is a {\it left Manin matrix\/};
see \cite[Prop.~18]{cfr:ap}.
\epf

An alternative (longer) way to proof Proposition~\ref{prop:hchan}
is to derive it from the character formulas of Proposition~\ref{prop:hcha}.
Indeed, $\di_u+E(u)$ coincides with the image of the matrix $T(u)\tss e^{\di_u}-1$
in the component of degree $-1$ of the graded algebra associated with the Yangian.
Here we extend the filtration on the Yangian to
the algebra of formal series $\Y(\gl_N)[[u^{-1},\di_u]]$ by setting
$\deg u^{-1}=\deg\di_u=-1$ so that the associated graded algebra is
isomorphic to $\U\big(\gl_N[t]\big)[[u^{-1},\di_u]]$.
Hence, for instance, the element on the left hand side of \eqref{hchan} can be found
as the image of the component of degree $-m$ of the expression
\ben
\tr\ts A^{(m)} \big(T_1(u)\tss e^{\di_u}-1\big)\dots
\big(T_m(u)\tss e^{\di_u}-1\big).
\een
The image of this expression under
the homomorphism
\eqref{yhch} can be found from \eqref{imbetheh}.

There is no known
analogue of the argument which we used in the proof of Proposition~\ref{prop:hchan}
for the $B$, $C$ and $D$ types.
Therefore to prove its counterparts for these types we have to resort
to the argument making use of the character formulas
of Sec.~\ref{subsec:yagn}.

\subsection{Types $B$, $C$ and $D$}
\label{subsec:csbcd}

Recall that $F_{ij}[r]=F_{ij}\ts t^r$ with $r\in\ZZ$ denote
elements of the loop algebra $\g_N[t,t^{-1}]$, where the $F_{ij}$
are standard generators of $\g_N$;
see Sec.~\ref{sec:cy}.

Consider the ascending filtration on the Yangian $\Y(\g_N)$ defined by
\ben
\deg t_{ij}^{(r)}=r-1.
\een
Denote by $\bar t_{ij}^{\ts(r)}$ the image of the generator
$t_{ij}^{(r)}$ in the $(r-1)$-th component of the
associated graded algebra $\gr\ts\Y(\g_N)$. By \cite[Theorem~3.6]{amr:rp}
the mapping
\ben
F_{ij}[r]\mapsto \bar t_{ij}^{\ts(r+1)},\qquad r\geqslant 0,
\een
defines an algebra isomorphism $\U\big(\g_N[t]\big)\to\gr\ts\Y(\g_N)$.
Our goal here is to use this isomorphism and Propositions~\ref{prop:hchb},
\ref{prop:hchd} and \ref{prop:hchc} to calculate the Harish-Chandra images
of certain elements of $\U\big(\g_N[t]\big)$ defined with the use of
the corresponding operators \eqref{symo} and \eqref{symsp}. These elements
generate a commutative subalgebra of $\U\big(\g_N[t]\big)$ and they
can be obtained from the generators \eqref{deftr} of the Feigin--Frenkel center
by an application of the vertex algebra structure on the vacuum module
$V_{-h^{\vee}}(\g_N)$; see \cite[Sec.~5]{m:ff}.

We identify the universal enveloping algebra $\U(\g_N)$
with a subalgebra of $\U\big(\g_N[t]\big)$ via the embedding
$F_{ij}\mapsto F_{ij}[0]$. Then $\U\big(\g_N[t]\big)$ can be regarded as a $\g_N$-module
with the adjoint action. Denote by $\U\big(\g_N[t]\big)^{\h}$ the
subalgebra of $\h$-invariants under this action.
Consider the left ideal $I$ of the algebra $\U\big(\g_N[t]\big)$ generated by
all elements $F_{ij}[r]$ with the conditions
$1\leqslant i<j\leqslant N$ and $r\geqslant 0$.
By the Poincar\'e--Birkhoff--Witt
theorem, the intersection
$\U\big(\g_N[t]\big)^{\h}\cap I$ is a two-sided ideal of
$\U\big(\g_N[t]\big)^{\h}$. Moreover,
the quotient of $\U\big(\g_N[t]\big)^{\h}$ by this ideal is isomorphic
to the commutative algebra freely generated by the images
of the elements $F_{ii}[r]$ with $i=1,\dots,n$ and $r\geqslant 0$
in the quotient. We will write $\mu_i[r]$
for this image of $F_{ii}[r]$ and extend this notation to
all values $i=1,\dots,N$ so that $\mu_{i'}[r]=-\mu_{i}[r]$
for all $i$. We get an analogue of the Harish-Chandra
homomorphism \eqref{hchcl},
\beql{yhchcu}
\U\big(\g_N[t]\big)^{\h}\to\CC[\mu_i[r]\ts|\ts i=1,\dots,n,\ r\geqslant 0].
\eeq
We will combine the elements $F_{ij}[r]$ into the formal series
\ben
F_{ij}(u)=\sum_{r=0}^{\infty} F_{ij}[r]\ts u^{-r-1}
\een
and write
\ben
\mu_i(u)=\sum_{r=0}^{\infty} \mu_i[r]\ts u^{-r-1},\qquad i=1,\dots,N.
\een
Then $\mu_i(u)$
is understood as the image of the series
$F_{ii}(u)$ under the homomorphism \eqref{yhchcu}.

It is clear from the definitions of the homomorphisms \eqref{yhchbcd} and
\eqref{yhchcu}, that the graded version of \eqref{yhchbcd}
coincides with \eqref{yhchcu} in the sense that the following diagram commutes
\beql{cdhch}
\begin{CD}
\U\big(\g_N[t]\big)^{\h} @>>> \CC\big[\mu_i[r]\big]\\
@VVV     @VVV\\
\gr\ts\Y(\g_N)^{\h} @>>> \gr\ts\CC\big[\la^{(r+1)}_i\big],
\end{CD}
\eeq
where $i$ ranges over the set $\{1,\dots,n\}$ while $r\geqslant 0$ and
the second vertical arrow indicates the isomorphism which takes
$\mu_i[r]$ to the image of $\la^{(r+1)}_i$ in the graded polynomial algebra
with the grading defined by the assignment $\deg \la^{(r+1)}_i=r$.

In what follows we extend the filtration on the Yangian to
the algebra of formal series $\Y(\g_N)[[u^{-1},\di_u]]$ by setting
$\deg u^{-1}=\deg\di_u=-1$. The associated graded algebra will then be
isomorphic to $\U\big(\g_N[t]\big)[[u^{-1},\di_u]]$.
We consider tensor product algebras
\beql{tenprkcu}
\underbrace{\End\CC^{N}\ot\dots\ot\End\CC^{N}}_m{}\ot \U\big(\g_N[t]\big)[[u^{-1},\di_u]]
\eeq
and use matrix notation as in \eqref{matnot}.

\subsubsection{Series $B_n$}

Take $\g_N=\oa_{N}$ with $N=2n+1$ and consider the operator
$S^{(m)}$ defined in \eqref{symo}. We also use notation
\eqref{galm} with $\om=N$ and \eqref{comp}.
The trace
is understood to be taken over all copies of the endomorphism
algebra $\End\CC^N$ in \eqref{tenprkcu}.

\bth\label{thm:hchbn}
For the image under the Harish-Chandra homomorphism \eqref{yhchcu} we have
\begin{multline}\label{hchbn}
\ga_m(N)\ts\tr\ts S^{(m)} \big(\di_u+F_1(u)\big)\dots \big(\di_u+F_m(u)\big)\\[1em]
{}\mapsto
h_m\big(\di_u+\mu_{1}(u),\dots,\di_u+\mu_{n}(u),
\di_u+\mu_{n'}(u),\dots,\di_u+\mu_{1'}(u)\big).
\end{multline}
\eth

\bpf
The element $\di_u+F(u)$ coincides with the image of the matrix $T(u)\tss e^{\di_u}-1$
in the component of degree $-1$ of the graded algebra associated with the Yangian.
Therefore the element on the left hand side of \eqref{hchbn} can be found
as the image of the component of degree $-m$ of the expression
\beql{yantb}
\ga_m(N)\ts\tr\ts S^{(m)} \big(T_1(u)\tss e^{\di_u}-1\big)\dots
\big(T_m(u)\tss e^{\di_u}-1\big).
\eeq
Hence, the theorem can be proved by making use of
the commutative diagram \eqref{cdhch} and
the Harish-Chandra image
of \eqref{yantb} implied by Proposition~\ref{prop:hchb}.
We have
\begin{multline}
\tr\ts S^{(m)}\ts \big(T_1(u)\tss e^{\di_u}-1\big)\dots
\big(T_m(u)\tss e^{\di_u}-1\big)\\
{}=\sum_{k=0}^m (-1)^{m-k}\ts \sum_{1\leqslant a_1<\dots<a_k\leqslant m}
\tr\ts S^{(m)}\ts T_{a_1}(u)\tss e^{\di_u}\dots T_{a_k}(u)\tss e^{\di_u}.
\non
\end{multline}
Each product $T_{a_1}(u)\tss e^{\di_u}\dots T_{a_k}(u)\tss e^{\di_u}$
can be written as $P\ts T_1(u)\tss e^{\di_u}\dots T_k(u)\tss e^{\di_u}\ts P^{-1}$,
where $P$ is the image in \eqref{tenprkcu} (with the identity component
in the last tensor factor)
of a permutation $p\in\Sym_m$
such that $p(r)=a_r$ for $r=1,\dots,k$. By the second property in
\eqref{sqpo} and the cyclic property of trace, we can bring
the above expression to the form
\ben
\sum_{k=0}^m (-1)^{m-k}\binom{m}{k} \ts
\tr\ts S^{(m)}\ts T_{1}(u)\tss e^{\di_u}\dots T_{k}(u)\tss e^{\di_u}.
\een
Now apply \cite[Lemma~4.1]{m:ff} to calculate
the partial traces of the symmetrizer $S^{(m)}$ over the copies
$k+1,\dots,m$ of the algebra $\End\CC^N$ in \eqref{tenprke} to get
\ben
\tr^{}_{k+1,\dots,m}\ts S^{(m)}
=\frac{\ga_k(N)}{\ga_m(N)}\ts \binom{N+m-2}{m-k}\ts\binom{m}{k}^{-1}\ts S^{(k)}.
\een
Thus, by Proposition~\ref{prop:hchb}, the Harish-Chandra image of the
expression \eqref{yantb} is found by
\beql{lahchb}
\sum_{k=0}^m(-1)^{m-k}\ts\ga_k(N)\ts\binom{N+m-2}{m-k}
\sum_{1\leqslant i_1\leqslant\dots\leqslant i_k\leqslant N}
\la_{i_1}(u)\tss e^{\di_u}\dots \la_{i_k}(u)\tss e^{\di_u}
\eeq
with the condition that $n+1$ occurs among the summation indices $i_1,\dots,i_k$
at most once.

The next step is to express \eqref{lahchb} in terms of the new variables
\beql{sila}
\si_i(u)=\la_{i}(u)\tss e^{\di_u}-1,\qquad i=1,\dots,N.
\eeq
This is done by a combinatorial argument as shown in the following lemma.

\ble\label{lem:lasib}
The expression \eqref{lahchb} multiplied by $-2\ts\binom{N/2-2}{N+m-2}$
equals
\begin{multline}
\sum_{r=0}^m\binom{N/2-2}{N+r-3}\sum_{a_1+\dots+a_{1'}=r}
\si_{1}(u)^{a_1}\dots \si_{n}(u)^{a_n}
\si_{n'}(u)^{a_{n'}}\dots \si_{1'}(u)^{a_{1'}}\\[1em]
{}+\sum_{r=1}^m\binom{N/2-2}{N+r-3}\sum_{a_1+\dots+a_{1'}=r-1}
\si_{1}(u)^{a_1}\dots \si_{n}(u)^{a_n}\big(\si_{n+1}(u)+2\big)\ts
\si_{n'}(u)^{a_{n'}}\dots \si_{1'}(u)^{a_{1'}},
\non
\end{multline}
where $a_1,\dots,a_{1'}$ run over nonnegative integers.
\ele

\bpf
The statement is verified by substituting \eqref{sila} into both terms
and calculating the coefficients of the sum
\beql{sumla}
\sum_{1\leqslant i_1\leqslant\dots\leqslant i_k\leqslant N}
\la_{i_1}(u)\tss e^{\di_u}\dots \la_{i_k}(u)\tss e^{\di_u}
\eeq
for all $0\leqslant k\leqslant m$, where $n+1$ occurs among
the summation indices $i_1,\dots,i_k$
at most once.
Note
the following expansion formula for the noncommutative complete symmetric
functions \eqref{comp},
\beql{comexp}
h_r(x_1-1,\dots,x_p-1)=\sum_{k=0}^r (-1)^{r-k}\binom{p+r-1}{r-k}\ts
h_k(x_1,\dots,x_p).
\eeq
Take $x_i=\la_i(u)\tss e^{\di_u}$ with $i=1,\dots,n,n'\dots,1'$
and apply \eqref{comexp} with $p=2n$ to the first term in the expression
of the lemma. Using a similar expansion for the second term
we find that the coefficient of the sum \eqref{sumla}
in the entire expression will be found as
\ben
\sum_{r=k}^m \binom{N/2-2}{N+r-3}\binom{N+r-3}{r-k}
=\binom{N/2-2}{N+k-3}\binom{N/2+m-1}{m-k},
\een
which coincides with
\ben
-2\ts (-1)^{m-k}\ts\ga_k(N)\ts\binom{N/2-2}{N+m-2}\ts\binom{N+m-2}{m-k},
\een
as claimed.
\epf

Denote the expression in Lemma~\ref{lem:lasib} by $A_m$.
Since the degree of
the element \eqref{yantb} is $-m$, its Harish-Chandra image \eqref{lahchb}
and the expression $A_m$ also have degree $-m$.
Observe that the terms in the both sums of $A_m$ are independent of $m$ so that
$A_{m+1}=A_m+B_{m+1}$, where
\begin{multline}
B_{m+1}=\binom{N/2-2}{N+m-2}\sum_{a_1+\dots+a_{1'}=m+1}
\si_{1}(u)^{a_1}\dots \si_{n}(u)^{a_n}
\si_{n'}(u)^{a_{n'}}\dots \si_{1'}(u)^{a_{1'}}\\[1em]
{}+\binom{N/2-2}{N+m-2}\sum_{a_1+\dots+a_{1'}=m}
\si_{1}(u)^{a_1}\dots \si_{n}(u)^{a_n}\big(\si_{n+1}(u)+2\big)\ts
\si_{n'}(u)^{a_{n'}}\dots \si_{1'}(u)^{a_{1'}}.
\non
\end{multline}
Since $A_{m+1}$ has degree $-m-1$, its
component of degree $-m$ is zero, and so the sum of the homogeneous components
of degree $-m$ of $A_m$ and $B_{m+1}$ is zero.
However,
each element $\si_i(u)$ has degree $-1$
with the top degree component equal to $\di_u+\mu_i(u)$.
This implies that the
component of $A_m$ of degree $-m$ equals the component of degree $-m$
of the term
\ben
-2\ts\binom{N/2-2}{N+m-2}\ts\sum_{a_1+\dots+a_{1'}=m}
\si_{1}(u)^{a_1}\dots \si_{n}(u)^{a_n}
\si_{n'}(u)^{a_{n'}}\dots \si_{1'}(u)^{a_{1'}}.
\een
Taking into account the constant factor used in Lemma~\ref{lem:lasib},
we can conclude that the component in question coincides with the
noncommutative complete symmetric function as given in \eqref{hchbn}.
\epf

As we have seen in the proof of the theorem,
all components of the expression in Lemma~\ref{lem:lasib}
of degrees exceeding $-m$ are equal to zero. Since the summands do not
depend on $m$, we derive the following corollary.

\bco\label{cor:infcob}
The series
\begin{multline}
\sum_{r=0}^{\infty}\binom{N/2-2}{N+r-3}\sum_{a_1+\dots+a_{1'}=r}
\si_{1}(u)^{a_1}\dots \si_{n}(u)^{a_n}
\si_{n'}(u)^{a_{n'}}\dots \si_{1'}(u)^{a_{1'}}\\[1em]
{}+\sum_{r=1}^{\infty}\binom{N/2-2}{N+r-3}\sum_{a_1+\dots+a_{1'}=r-1}
\si_{1}(u)^{a_1}\dots \si_{n}(u)^{a_n}\big(\si_{n+1}(u)+2\big)\ts
\si_{n'}(u)^{a_{n'}}\dots \si_{1'}(u)^{a_{1'}}
\non
\end{multline}
is equal to zero.
\qed
\eco

\subsubsection{Series $D_n$}

Now take $\g_N=\oa_{N}$ with $N=2n$ and consider the operator
$S^{(m)}$ defined in \eqref{symo}. We keep using notation
\eqref{galm} with $\om=N$ and \eqref{comp}.

\bth\label{thm:hchdn}
For the image under the Harish-Chandra homomorphism \eqref{yhchcu} we have
\begin{multline}
\ga_m(N)\ts\tr\ts S^{(m)} \big(\di_u+F_1(u)\big)\dots \big(\di_u+F_m(u)\big)\\[1em]
{}\mapsto
{\textstyle\frac12}\ts h_m\big(\di_u+\mu_{1}(u),\dots,\di_u+\mu_{n-1}(u),
\di_u+\mu_{n'}(u),\dots,\di_u+\mu_{1'}(u)\big)\\[0.5em]
{}+{\textstyle\frac12}\ts h_m\big(\di_u+\mu_{1}(u),\dots,\di_u+\mu_{n}(u),
\di_u+\mu_{(n-1)'}(u),\dots,\di_u+\mu_{1'}(u)\big).
\non
\end{multline}
\eth

\bpf
We repeat the beginning of the proof of Theorem~\ref{thm:hchdn}
with $N$ now taking the even value $2n$ up to the application
of the formula for Yangian characters. This time
we apply Proposition~\ref{prop:hchd} to conclude that the Harish-Chandra image of the
expression \eqref{yantb} is found by
\beql{lahchd}
\sum_{k=0}^m(-1)^{m-k}\ts\ga_k(2n)\ts\binom{2n+m-2}{m-k}
\sum_{1\leqslant i_1\leqslant\dots\leqslant i_k\leqslant 2n}
\la_{i_1}(u)\tss e^{\di_u}\dots \la_{i_k}(u)\tss e^{\di_u}
\eeq
with the condition that $n$ and $n'$ do not occur simultaneously
among the summation indices $i_1,\dots,i_k$.
Introducing new variables by the same formulas \eqref{sila}
we come to the $D_n$ series counterpart of Lemma~\ref{lem:lasib},
where we use the notation
\ben
c_r=(-1)^{r-1}\ts\binom{2\tss n+r-2}{n-1}^{-1}.
\een

\ble\label{lem:lasid}
The expression \eqref{lahchd} multiplied by $2\tss c_m$ equals
\begin{multline}
2\tss c_m\ts
\sum_{\underset{\scriptstyle a^{}_n=a_{n'}=0}{a_1+\dots+a_{1'}=m}}
\si_{1}(u)^{a_1}\dots \si_{1'}(u)^{a_{1'}}
+c_m\ts
\sum_{\underset{\scriptstyle \text{only one of\ } a^{}_n\text{\ and\ }a_{n'}
\text{\ is zero}}{a_1+\dots+a_{1'}=m}}
\si_{1}(u)^{a_1}\dots \si_{1'}(u)^{a_{1'}}\\[0.5em]
{}-\sum_{r=1}^m\ts\frac{r\ts c_r}{n+r-1}\ts
\sum_{\underset{\scriptstyle a^{}_n=a_{n'}=0}{a_1+\dots+a_{1'}=r}}
\si_{1}(u)^{a_1}\dots \si_{1'}(u)^{a_{1'}}\\[0.5em]
{}+\sum_{r=1}^m\ts\frac{(n-1)\ts c_r}{n+r-1}\ts
\sum_{\underset{\scriptstyle \text{only one of\ } a^{}_n\text{\ and\ }a_{n'}
\text{\ is zero}}{a_1+\dots+a_{1'}=r}}
\si_{1}(u)^{a_1}\dots \si_{1'}(u)^{a_{1'}},
\non
\end{multline}
where $a_1,\dots,a_{1'}$ run over nonnegative integers.
\ele

\bpf
Substitute \eqref{sila} into the expression
and calculate the coefficients of the sum
\beql{sumlad}
\sum_{1\leqslant i_1\leqslant\dots\leqslant i_k\leqslant 2n}
\la_{i_1}(u)\tss e^{\di_u}\dots \la_{i_k}(u)\tss e^{\di_u}.
\eeq
The argument splits into two cases, depending on whether
neither of $n$ and $n'$ occurs among the summation indices $i_1,\dots,i_k$
in \eqref{sumlad} or only one of them occurs. The application
of the expansion formula \eqref{comexp} brings this to a straightforward
calculation with the binomial coefficients in both cases.
\epf

Let $A_m$ denote the four-term expression in Lemma~\ref{lem:lasid}. This expression
equals $2\tss c_m$ times the Harish-Chandra image of \eqref{yantb} and so
$A_m$ has degree $-m$. Hence, the component of degree $-m$
of the expression $A_{m+1}$ is zero.
On the other hand,
each element $\si_i(u)$ has degree $-1$
with the top degree component equal to $\di_u+\mu_i(u)$.
This implies that the component of degree $-m$ in the sum
of the third and fourth terms in $A_m$ is zero. Therefore, the
component of $A_m$ of degree $-m$ equals the component of degree $-m$
in the sum of the first and the second terms.
Taking into account the constant factor $2\tss c_m$,
we conclude that the component takes the desired form.
\epf

The following corollary is implied by the proof of the theorem.

\bco\label{cor:infcod}
The series
\begin{multline}
{}-\sum_{r=1}^{\infty}\ts\frac{r\ts c_r}{n+r-1}\ts
\sum_{\underset{\scriptstyle a^{}_n=a_{n'}=0}{a_1+\dots+a_{1'}=r}}
\si_{1}(u)^{a_1}\dots \si_{1'}(u)^{a_{1'}}\\[0.5em]
{}+\sum_{r=1}^{\infty}\ts\frac{(n-1)\ts c_r}{n+r-1}\ts
\sum_{\underset{\scriptstyle \text{only one of\ } a^{}_n\text{\ and\ }a_{n'}
\text{\ is zero}}{a_1+\dots+a_{1'}=r}}
\si_{1}(u)^{a_1}\dots \si_{1'}(u)^{a_{1'}}
\non
\end{multline}
is equal to zero.
\qed
\eco

\subsubsection{Series $C_n$}

Now we let $\g_N=\spa_{N}$ with $N=2n$ and consider the operator
$S^{(m)}$ defined in \eqref{symsp}. We also use notation
\eqref{galm} with $\om=-2n$ and \eqref{elem}.
Although the operator $S^{(m)}$ is defined only for
$m\leqslant n+1$, it is possible to extend the values of expressions
of the form \eqref{deftr} and those which are used in the next theorem
to all $m$ with $m\leqslant 2n+1$; see
\cite[Sec.~3.3]{m:ff}. The Harish-Chandra images turn out to be given by
the same expression for all these values of $m$.
We postpone the proof to Corollary~\ref{cor:extva} below, and assume
first that $m\leqslant n$.

\bth\label{thm:hchcn}
For all $1\leqslant m\leqslant n$
for the image under the Harish-Chandra homomorphism \eqref{yhchcu} we have
\begin{multline}\label{hchcn}
\ga_m(-2n)\ts\tr\ts S^{(m)} \big(\di_u-F_1(u)\big)\dots \big(\di_u-F_m(u)\big)\\[1em]
{}\mapsto
e_m\big(\di_u+\mu_{1}(u),\dots,\di_u+\mu_{n}(u),\di_u,
\di_u+\mu_{n'}(u),\dots,\di_u+\mu_{1'}(u)\big).
\end{multline}
\eth

\bpf
The element $\di_u-F(u)$ coincides with the image of the matrix $1-T(u)\tss e^{-\di_u}$
in the component of degree $-1$ of the graded algebra associated with the Yangian.
Hence the left hand side of \eqref{hchcn} can be found
as the image of the component of degree $-m$ of the expression
\beql{yantc}
(-1)^m\tss\ga_m(-2\tss n)\ts\tr\ts S^{(m)} \big(T_1(u)\tss e^{-\di_u}-1\big)\dots
\big(T_m(u)\tss e^{-\di_u}-1\big).
\eeq
Now we use
the commutative diagram \eqref{cdhch} and
the Harish-Chandra image
of \eqref{yantc} implied by Proposition~\ref{prop:hchc}.
We have
\begin{multline}
\tr\ts S^{(m)}\ts \big(T_1(u)\tss e^{-\di_u}-1\big)\dots
\big(T_m(u)\tss e^{-\di_u}-1\big)\\
{}=\sum_{k=0}^m (-1)^{m-k}\ts \sum_{1\leqslant a_1<\dots<a_k\leqslant m}
\tr\ts S^{(m)}\ts T_{a_1}(u)\tss e^{-\di_u}\dots T_{a_k}(u)\tss e^{-\di_u}.
\label{wga}
\end{multline}
As in the proof of Theorem~\ref{thm:hchcn}, we use
the second property in
\eqref{sqpo} and the cyclic property of trace to bring
the expression to the form
\ben
\sum_{k=0}^m (-1)^{m-k}\binom{m}{k} \ts
\tr\ts S^{(m)}\ts T_{1}(u)\tss e^{-\di_u}\dots T_{k}(u)\tss e^{-\di_u}.
\een
Further, the partial traces of the symmetrizer $S^{(m)}$ over the copies
$k+1,\dots,m$ of the algebra $\End\CC^N$ in \eqref{tenprke}
are found by applying \cite[Lemma~4.1]{m:ff} to get
\ben
\tr^{}_{k+1,\dots,m}\ts S^{(m)}
=\frac{\ga_k(-2\tss n)}{\ga_m(-2\tss n)}\ts
\binom{2\tss n-k+1}{m-k}\ts\binom{m}{k}^{-1}\ts S^{(k)}.
\een
By Proposition~\ref{prop:hchc} and Corollary~\ref{cor:eqka},
the Harish-Chandra image of the
expression \eqref{yantc} is found by
\beql{lahchc}
\sum_{k=0}^m(-1)^{k}\ts\ga_k(-2\tss n)\ts\binom{2\tss n-k+1}{m-k}\ts
\sum_{1\leqslant i_1<\dots< i_k\leqslant 2n+2}
\vk_{i_1}(u)\tss e^{-\di}\dots \vk_{i_k}(u)\tss e^{-\di}.
\eeq
Introduce new variables by
\beql{silac}
\si_i(u)=\vk_i(u)\tss e^{-\di}-1,\qquad i=1,\dots,2n+2,\quad i\ne n+2,
\eeq
and $\si_{n+2}(u)=\vk_{n+2}(u)\tss e^{-\di}+1$.

\ble\label{lem:lasic}
For $m\leqslant n$
the expression \eqref{lahchc} multiplied by $2\tss(-1)^m\tss\binom{2\tss n-m+1}{n+1}$
equals
\begin{multline}
\sum_{r=0}^m\binom{2\tss n-r+2}{n+1}\sum_{1\leqslant i_1<\dots< i_r\leqslant 2n+2}
\si_{i_1}(u)\dots \si_{i_r}(u)\\[0.5em]
{}-2\ts \sum_{r=0}^{m-1}\binom{2\tss n-r+1}{n+1}
\sum_{\underset{\scriptstyle i_s\ne n+2}
{1\leqslant i_1<\dots< i_r\leqslant 2n+2}}
\si_{i_1}(u)\dots \si_{i_r}(u),
\non
\end{multline}
where $n+2$ does not occur among the summation indices in the last sum.
\ele

\bpf
Substituting \eqref{silac} into the expression and simplifying gives
\beql{sumvk}
\sum_{r=0}^m\binom{2\tss n-r+2}{n+1}\sum_{1\leqslant i_1<\dots< i_r\leqslant 2n+2}
\big(\vk_{i_1}(u)\tss e^{-\di}-1\big)\dots
\big(\vk_{i_r}(u)\tss e^{-\di}-1\big).
\eeq
Now use the expansion formula for the noncommutative elementary symmetric
functions \eqref{elem},
\ben
e_r(x_1-1,\dots,x_p-1)=\sum_{k=0}^r (-1)^{r-k}\binom{p-k}{r-k}\ts
e_k(x_1,\dots,x_p).
\een
Taking $x_i=\vk_i(u)\tss e^{-\di_u}$ with $i=1,\dots,2n+2$,
it is straightforward to verify that
the coefficient of the sum
\ben
\sum_{1\leqslant i_1<\dots< i_k\leqslant N}
\vk_{i_1}(u)\tss e^{-\di_u}\dots \vk_{i_k}(u)\tss e^{-\di_u}
\een
in \eqref{sumvk} equals
\ben
(-1)^{m-k}\binom{n-k}{m-k}\binom{2\tss n-k+2}{n+1}
\een
which coincides with
\ben
2\tss(-1)^{m-k}\tss\ga_k(-2\tss n)\tss
\binom{2\tss n-m+1}{n+1}\ts\binom{2\tss n-k+1}{m-k}
\een
as claimed.
\epf

For $m\leqslant n$ let $A_m$ denote the expression in Lemma~\ref{lem:lasic}.
Note that $A_m$ coincides with the Harish-Chandra image of \eqref{wga}
multiplied by $\binom{2\tss n-m+2}{n+1}$.
The proof of Lemma~\ref{lem:lasic} and the second part of Corollary~\ref{cor:eqka}
show that $A_m$ is also well-defined for the value $m=n+1$ and $A_{n+1}=0$.

Since the degree of
the element \eqref{yantc} is $-m$, for $m\leqslant n$
the expression $A_m$ also has degree $-m$. Hence, the component of degree $-m$
of the expression $A_{m+1}$ is zero; this holds for $m=n$ as well, because
$A_{n+1}=0$. Furthermore, each element $\si_i(u)$ has degree $-1$ and so
the component of $A_m$ of degree $-m$ must be equal to the component of degree $-m$
of the expression
\ben
2\ts \binom{2\tss n-m+1}{n+1}
\sum_{\underset{\scriptstyle i_s\ne n+2}
{1\leqslant i_1<\dots< i_m\leqslant 2n+2}}
\si_{i_1}(u)\dots \si_{i_m}(u).
\een
The component of $\si_i(u)$ of degree $-1$ equals
\ben
\begin{cases}-\di_u+\mu_i(u)\qquad&\text{for}\quad i=1,\dots,n,\\
-\di_u\qquad&\text{for}\quad i=n+1,\\
-\di_u+\mu_{i-2}(u)\qquad&\text{for}\quad i=n+3,\dots,2n+2.
\end{cases}
\een
The proof is completed by taking the signs and the constant
factor used in Lemma~\ref{lem:lasic} into account.
\epf

\section{Classical $\Wc$-algebras}
\label{sec:cw}
\setcounter{equation}{0}

We define the classical $\Wc$-algebra
$\Wc(\g)$ associated with a simple Lie algebra $\g$
following \cite[Sec.~8.1]{f:lc},
where more details
and proofs can be found.
We let $\h$ denote a Cartan subalgebra of $\g$ and let $\mu_1,\dots,\mu_n$
be a basis of $\h$. The universal enveloping algebra $\U(t^{-1}\h[t^{-1}])$
will be identified with the algebra of polynomials in the infinitely
many variables $\mu_i[r]$ with $i=1,\dots,n$ and $r<0$
and will be denoted by $\pi_0$. We will also use the extended
algebra with the additional generator $\tau$ subject to the relations
\ben
\big[\tau,\mu_i[r]\tss\big]=-r\ts \mu_i[r-1],
\een
implied by \eqref{taur}. The extended algebra
is isomorphic to $\pi_0\ot\CC[\tau]$ as a vector space.
Furthermore, we will need the operator
$T={\rm ad}\ts\tau$ which is the derivation $T:\pi_0\to\pi_0$
defined on the generators by the relations
\ben
T\ts\mu_i[r]=-r\ts \mu_i[r-1].
\een
In particular, $T\ts 1=0$.
The {\it classical $\Wc$-algebra\/}
is defined as the subspace $\Wc(\g)\subset\pi_0$ spanned by
the elements which are annihilated
by the {\it screening operators\/}
\ben
V_i:\pi_0\to \pi_0, \qquad i=1,\dots,n,
\een
which we will write down explicitly for each classical
type below,\footnote{Our $V_i$ essentially coincides
with the operator $\overline V_i[1]$ in the notation of \cite[Sec.~7.3.4]{f:lc},
which is associated with the Langlands dual Lie algebra ${}^L\g$.}
\ben
\Wc(\g)=\{P\in \pi_0\ |\ V_i\ts P=0,\quad i=1,\dots,n\}.
\een
The operators $V_i$ are derivations of $\pi_0$ so that $\Wc(\g)$
is a subalgebra of $\pi_0$.
The subalgebra $\Wc(\g)$ is $T$-invariant. Moreover, there exist
elements $B_1,\dots,B_n\in \Wc(\g)$ such that the family of elements
$T^{\tss r}B_i$ with $i=1,\dots,n$ and $r\geqslant 0$ is algebraically
independent and generates the algebra $\Wc(\g)$. We will call
$B_1,\dots,B_n$ a {\it complete set of generators\/} of $\Wc(\g)$.
Examples of such sets in the classical types will be given below.

We extend the screening operators to the algebra $\pi_0\ot\CC[\tau]$
by
\ben
V_i\big(P\ot Q(\tau)\big)=V_i(P)\ot Q(\tau),\qquad P\in\pi_0,\quad Q(\tau)\in\CC[\tau].
\een

\subsection{Screening operators and generators for $\Wc(\gl_N)$}
\label{subsec:wagln}

Here $\pi_0$ is the algebra of polynomials in the variables $\mu_i[r]$
with $i=1,\dots,N$ and $r<0$. The screening operators $V_1,\dots,V_{N-1}$
are defined by
\ben
V_i=\sum_{r=0}^{\infty} V_{i\ts [r]}\ts
\Big(\frac{\di}{\di\mu_i[-r-1]}-\frac{\di}{\di\mu_{i+1}[-r-1]}\Big),
\een
where the coefficients $V_{i\ts [r]}$ are found from the expansion
of a formal generating function in a variable $z$,
\ben
\sum_{r=0}^{\infty} V_{i\ts [r]}\ts z^r=\exp\ts\sum_{m=1}^{\infty}
\frac{\mu_i[-m]-\mu_{i+1}[-m]}{m}\ts z^m.
\een
Define elements $\Ec_1,\dots,\Ec_N$ of $\pi_0$ by
the expansion in $\pi_0\ot\CC[\tau]$,
\beql{epab}
\big(\tau+\mu_N[-1]\big)\dots\big(\tau+\mu_1[-1]\big)
=\tau^N+ \Ec_1\ts\tau^{N-1}+\dots+\Ec_N,
\eeq
known as the {\it Miura transformation\/}.
Explicitly, using the notation \eqref{elem}
we can write the coefficients as
\beql{bmgln}
\Ec_m=e_m\big(T+\mu_1[-1],\dots,T+\mu_N[-1]\big),
\eeq
which follows easily from \eqref{epab} by induction.
The family $\Ec_1,\dots,\Ec_N$ is a complete set of generators of $\Wc(\gl_N)$.
Verifying that all elements $\Ec_i$ are annihilated by the screening
operators is straightforward. This is implied by the relations
for the operators on $\pi_0$,
\beql{vit}
V_i\ts T=\big(T+\mu_i[-1]-\mu_{i+1}[-1]\big)\tss V_i,\qquad i=1,\dots,N-1.
\eeq
They imply the corresponding relations for the operators on $\pi_0\ot\CC[\tau]$,
\beql{vitau}
V_i\ts \tau=\big(\tau+\mu_i[-1]-\mu_{i+1}[-1]\big)\tss V_i,\qquad i=1,\dots,N-1,
\eeq
where $\tau$ is regarded as the operator of left multiplication by $\tau$.
For each $i$ the relation
\ben
V_i\ts \big(\tau+\mu_N[-1]\big)\dots\big(\tau+\mu_1[-1]\big)=0
\een
then follows easily. Indeed, it reduces to the particular case $N=2$
where we have
\ben
\bal
V_1\ts \big(\tau+\mu_2[-1]\big)\big(\tau+\mu_1[-1]\big)
{}&=\Big(\big(\tau+\mu_1[-1]-\mu_2[-1]\big)\ts V_1+\mu_2[-1]\ts V_1-1\Big)
\big(\tau+\mu_1[-1]\big)\\
{}&=\big(\tau+\mu_1[-1]\big)\ts V_1\ts\big(\tau+\mu_1[-1]\big)-
\big(\tau+\mu_1[-1]\big)=0.
\eal
\een
Showing that the elements $T^{\tss r}\Ec_i$ are algebraically
independent generators requires a comparison of the sizes of
graded components of $\pi_0$ and $\Wc(\gl_N)$.

By the definitions \eqref{comp} and \eqref{elem}, we have
the relations
\beql{elemcomp}
\sum_{k=0}^m (-1)^k\ts\Ec_k\ts h_{m-k}\big(T+\mu_1[-1],\dots,T+\mu_N[-1]\big)=0
\eeq
for $m\geqslant 1$, where $\Ec_0=1$ and $\Ec_k=0$ for $k>N$. They imply
that all elements
\beql{hmgln}
h_m\big(T+\mu_1[-1],\dots,T+\mu_N[-1]\big),\qquad m\geqslant 1,
\eeq
belong to $\Wc(\gl_N)$. Moreover, the family \eqref{hmgln} with $m=1,\dots,N$
is a complete set of generators of $\Wc(\gl_N)$.

Note that the classical $\Wc$-algebra $\Wc(\sll_N)$ associated with
the special linear Lie algebra $\sll_N$ can be obtained as the quotient
of $\Wc(\gl_N)$ by the relation $\Ec_1=0$.

\subsection{Screening operators
 and generators for $\Wc(\oa_N)$ and $\Wc(\spa_N)$}
\label{subsec:wabcd}

Now $\pi_0$ is the algebra of polynomials in the variables $\mu_i[r]$
with $i=1,\dots,n$ and $r<0$. The families of generators
of the algebras $\Wc(\oa_N)$ and $\Wc(\spa_N)$ reproduced below
were constructed in \cite[Sec.~8]{ds:la},
where equations of the KdV type were introduced for arbitrary simple Lie algebras.
The generators are associated with
the Miura transformations of the corresponding equations.

\subsubsection{Series $B_n$}

The screening operators $V_1,\dots,V_n$
are defined by
\beql{vigen}
V_i=\sum_{r=0}^{\infty} V_{i\ts [r]}\ts
\Big(\frac{\di}{\di\mu_i[-r-1]}-\frac{\di}{\di\mu_{i+1}[-r-1]}\Big),
\eeq
for $i=1,\dots,n-1$, and
\ben
V_n=\sum_{r=0}^{\infty} V_{n\ts [r]}\ts
\frac{\di}{\di\mu_n[-r-1]},
\een
where the coefficients $V_{i\ts [r]}$ are found from the expansions
\ben
\sum_{r=0}^{\infty} V_{i\ts [r]}\ts z^r=\exp\ts\sum_{m=1}^{\infty}
\frac{\mu_i[-m]-\mu_{i+1}[-m]}{m}\ts z^m,\qquad i=1,\dots,n-1
\een
and
\ben
\sum_{r=0}^{\infty} V_{n\ts [r]}\ts z^r=\exp\ts\sum_{m=1}^{\infty}
\frac{\mu_n[-m]}{m}\ts z^m.
\een
Define elements $\Ec_2,\dots,\Ec_{2n+1}$ of $\pi_0$ by
the expansion
\begin{multline}\label{epabb}
\big(\tau-\mu_1[-1]\big)\dots\big(\tau-\mu_n[-1]\big)\ts\tau\ts\big(\tau+\mu_n[-1]\big)
\dots\big(\tau+\mu_1[-1]\big)\\[0.5em]
{}=\tau^{2n+1}+ \Ec_2\ts\tau^{2n-1}+\Ec_3\ts\tau^{2n-2}+\dots+\Ec_{2n+1}.
\end{multline}
All of them belong to $\Wc(\oa_{2n+1})$. By \eqref{bmgln} we have
\beql{bmglnb}
\Ec_m=e_m\big(T+\mu_1[-1],\dots,T+\mu_n[-1],T,T-\mu_n[-1],\dots,T-\mu_1[-1]\big).
\eeq
The family $\Ec_2,\Ec_4,\dots,\Ec_{2n}$ is a complete set of generators of $\Wc(\oa_{2n+1})$.
The relation
\beql{vib}
V_i\ts \big(\tau-\mu_1[-1]\big)\dots\big(\tau-\mu_n[-1]\big)
\ts\tau\ts\big(\tau+\mu_n[-1]\big)
\dots\big(\tau+\mu_1[-1]\big)=0
\eeq
is verified for $i=1,\dots,n-1$ in the same way as for $\gl_N$
with the use of \eqref{vitau}. Furthermore,
\ben
V_n\ts \tau=\big(\tau+\mu_n[-1]\big)\tss V_n,
\een
so that
\ben
\bal
V_n\ts \big(\tau-\mu_n[-1]\big)\ts\tau\ts\big(\tau+\mu_n[-1]\big)
{}&=\big(\tau\tss V_n-1\big)\ts\tau\ts\big(\tau+\mu_n[-1]\big)\\
{}&=\tau\ts\big(\tau+\mu_n[-1]\big)\ts\big(\tau+2\tss\mu_n[-1]\big)\tss V_n,
\eal
\een
which implies that \eqref{vib} holds for $i=n$ as well.

By \eqref{elemcomp} all elements
\beql{hmb}
h_m\big(T+\mu_1[-1],\dots,T+\mu_n[-1],T,T-\mu_n[-1],\dots,T-\mu_1[-1]\big)
\eeq
belong to $\Wc(\oa_{2n+1})$. The family of elements \eqref{hmb}
with $m=2,4,\dots,2n$ forms another complete set of
generators of $\Wc(\oa_{2n+1})$.

\subsubsection{Series $C_n$}

The screening operators $V_1,\dots,V_n$
are defined by \eqref{vigen}
for $i=1,\dots,n-1$, and
\ben
V_n=\sum_{r=0}^{\infty} V_{n\ts [r]}\ts
\frac{\di}{\di\mu_n[-r-1]},
\een
where
\ben
\sum_{r=0}^{\infty} V_{n\ts [r]}\ts z^r=\exp\ts\sum_{m=1}^{\infty}
\frac{2\tss\mu_n[-m]}{m}\ts z^m.
\een
Define elements $\Ec_2,\dots,\Ec_{2n}$ of $\pi_0$ by
the expansion
\begin{multline}
\big(\tau-\mu_1[-1]\big)\dots\big(\tau-\mu_n[-1]\big)\big(\tau+\mu_n[-1]\big)
\dots\big(\tau+\mu_1[-1]\big)\\[0.5em]
{}=\tau^{2n}+ \Ec_2\ts\tau^{2n-2}+\Ec_3\ts\tau^{2n-3}+\dots+\Ec_{2n}.
\non
\end{multline}
All of them belong to $\Wc(\spa_{2n})$.
By \eqref{bmgln} we have
\ben
\Ec_m=e_m\big(T+\mu_1[-1],\dots,T+\mu_n[-1],T-\mu_n[-1],\dots,T-\mu_1[-1]\big).
\een
The family $\Ec_2,\Ec_4,\dots,\Ec_{2n}$ is a complete set of generators of $\Wc(\spa_{2n})$.
The relation
\beql{vic}
V_i\ts \big(\tau-\mu_1[-1]\big)\dots\big(\tau-\mu_n[-1]\big)\big(\tau+\mu_n[-1]\big)
\dots\big(\tau+\mu_1[-1]\big)=0
\eeq
is verified for $i=1,\dots,n-1$ in the same way as for $\gl_N$
with the use of \eqref{vitau}. In the case $i=n$ we have
\ben
V_n\ts \tau=\big(\tau+2\tss\mu_n[-1]\big)\tss V_n,
\een
so that
\ben
\bal
V_n\ts \big(\tau-\mu_n[-1]\big)\big(\tau+\mu_n[-1]\big)
{}&=\Big(\big(\tau+\mu_n[-1]\big)\tss V_n-1\Big)\big(\tau+\mu_n[-1]\big)\\
{}&=\big(\tau+\mu_n[-1]\big)\big(\tau+3\tss\mu_n[-1]\big)\tss V_n,
\eal
\een
and \eqref{vic} with $i=n$ also follows.

It follows from \eqref{elemcomp} that the elements
\ben
h_m\big(T+\mu_1[-1],\dots,T+\mu_n[-1],T-\mu_n[-1],\dots,T-\mu_1[-1]\big)
\een
with $m=2,4,\dots,2n$ form another complete set of
generators of $\Wc(\spa_{2n})$.

\subsubsection{Series $D_n$}
\label{subsubsec:dn}

The screening operators $V_1,\dots,V_n$
are defined by \eqref{vigen}
for $i=1,\dots,n-1$, and
\ben
V_n=\sum_{r=0}^{\infty} V_{n\ts [r]}\ts
\Big(\frac{\di}{\di\mu_{n-1}[-r-1]}+\frac{\di}{\di\mu_n[-r-1]}\Big)
\een
where
\ben
\sum_{r=0}^{\infty} V_{n\ts [r]}\ts z^r=\exp\ts\sum_{m=1}^{\infty}
\frac{\mu_{n-1}[-m]+\mu_n[-m]}{m}\ts z^m.
\een
Define elements $\Ec_2,\Ec_3,\dots$ of $\pi_0$ by
the expansion of the {\it pseudo-differential operator\/}
\begin{multline}
\big(\tau-\mu_1[-1]\big)\dots\big(\tau-\mu_n[-1]\big)\ts\tau^{-1}
\ts\big(\tau+\mu_n[-1]\big)
\dots\big(\tau+\mu_1[-1]\big)\\[0.5em]
{}=\tau^{2n-1}+ \sum_{k=2}^{\infty}\Ec_k\ts\tau^{2n-k-1}.
\non
\end{multline}
The coefficients $\Ec_k$ are calculated with the use of the relations
\ben
\tau^{-1}\mu_i[-r-1]=\sum_{k=0}^{\infty}\frac{(-1)^k\tss (r+k)!}{r!}
\ts \mu_i[-r-k-1]\tss \tau^{-k-1}.
\een
All the elements $\Ec_k$ belong to $\Wc(\oa_{2n})$.
Moreover, define $\Ec^{\tss\prime}_n\in\pi_0$ by
\beql{pfaffim}
\Ec^{\tss\prime}_n=
\big(\mu_{1}[-1]-T\big)\dots \big(\mu_{n}[-1]-T\big),
\eeq
so that this element coincides with \eqref{pfaim}.
The family $\Ec_2,\Ec_4,\dots,\Ec_{2n-2},\Ec^{\tss\prime}_n$
is a complete set of generators of $\Wc(\oa_{2n})$.
The identity
\beql{vid}
V_i\ts \big(\tau-\mu_1[-1]\big)\dots\big(\tau-\mu_n[-1]\big)
\ts\tau^{-1}\ts\big(\tau+\mu_n[-1]\big)
\dots\big(\tau+\mu_1[-1]\big)=0
\eeq
is verified
with the use of \eqref{vitau} and the additional relations
\ben
V_i\ts \tau^{-1}=\big(\tau+\mu_i[-1]-\mu_{i+1}[-1]\big)^{-1}\tss V_i,\qquad i=1,\dots,n-1,
\een
and
\beql{vntauin}
V_n\ts \tau^{-1}=\big(\tau+\mu_{n-1}[-1]+\mu_n[-1]\big)^{-1}\tss V_n.
\eeq
In comparison with the types $B_n$ and $C_n$, an additional
calculation is needed for the case $i=n$ in \eqref{vid}.
It suffices to take $n=2$.
We have
\ben
\bal
V_2\ts \big(\tau-\mu_1[-1]\big)&\big(\tau-\mu_2[-1]\big)
\ts\tau^{-1}\ts\big(\tau+\mu_2[-1]\big)\big(\tau+\mu_1[-1]\big)\\
{}&=\Big(\big(\tau+\mu_2[-1]\big)\tss V_2-1\Big)\big(\tau-\mu_2[-1]\big)
\ts\tau^{-1}\ts\big(\tau+\mu_2[-1]\big)\big(\tau+\mu_1[-1]\big)\\
{}&=\Big(\big(\tau+\mu_2[-1]\big)\big(\tau+\mu_1[-1]\big)\tss V_2-2\tss\tau\Big)
\ts\tau^{-1}\ts\big(\tau+\mu_2[-1]\big)\big(\tau+\mu_1[-1]\big).
\eal
\een
Furthermore, applying the operator $V_2$ we find
\ben
\bal
V_2\ts \big(\tau+\mu_2[-1]\big)\big(\tau+\mu_1[-1]\big)
{}&=\Big(\big(\tau+\mu_1[-1]+2\tss\mu_2[-1]\big)\tss V_2+1\Big)\big(\tau+\mu_1[-1]\big)\\
{}&=2\tss\big(\tau+\mu_1[-1]+\mu_2[-1]\big)
\eal
\een
and so by \eqref{vntauin},
\ben
V_2
\ts\tau^{-1}\ts\big(\tau+\mu_2[-1]\big)\big(\tau+\mu_1[-1]\big)
=2
\een
thus completing the calculation.

The relations
\ben
V_i\ts \big(\mu_{1}[-1]-T\big)\dots \big(\mu_{n}[-1]-T\big)=0,\qquad i=1,\dots,n,
\een
are verified with the use of \eqref{vit}.

\section{Generators of the $\Wc$-algebras}
\label{sec:gw}
\setcounter{equation}{0}

Here we prove the Main Theorem stated in the Introduction by deriving it from
Theorems~\ref{thm:hchbn}, \ref{thm:hchdn} and \ref{thm:hchcn}.

Choose a basis $X_1,\dots,X_d$ of the simple Lie algebra $\g$ and write
the commutation relations
\ben
[X_i,X_j]=\sum_{k=1}^d c_{ij}^{\ts k}\ts X_k
\een
with structure constants $c_{ij}^{\ts k}$.
Consider the Lie algebras $\g[t]$ and $t^{-1}\g[t^{-1}]$ and combine their generators
into formal series in $u^{-1}$ and $u$,
\ben
X_i(u)=\sum_{r=0}^{\infty} X_i[r]\ts u^{-r-1}\Fand
X_i(u)_+=\sum_{r=0}^{\infty} X_i[-r-1]\ts u^{r}.
\een
The commutation relations of these Lie algebras written in terms
of the formal series take the form
\begin{align}
(u-v)\ts [X_i(u),X_j(v)]&=-\sum_{k=1}^d c_{ij}^{\ts k}\ts
\big(X_k(u)-X_k(v)\big),
\non\\
(u-v)\ts [X_i(u)_+,X_j(v)_+]&=\sum_{k=1}^d c_{ij}^{\ts k}\ts
\big(X_k(u)_+-X_k(v)_+\big).
\non
\end{align}
Observe that the second family of commutation relations is obtained from the first
by replacing $X_i(u)$ with the respective series $-X_i(u)_+$.

On the other hand, in the classical types,
the elements of the universal enveloping algebra
$\U(\g[t])$ and their Harish-Chandra images calculated in
Proposition~\ref{prop:hchan} and
Theorems~\ref{thm:hchbn}, \ref{thm:hchdn} and \ref{thm:hchcn}
are all expressed in terms of the series of the form $X_i(u)$.
Therefore, the corresponding Harish-Chandra images of the elements
of the universal enveloping algebra $\U(t^{-1}\g[t^{-1}])$
are readily found from those theorems by replacing $X_i(u)$
with the respective series $-X_i(u)_+$.

To be consistent with the definition for the Wakimoto modules in
\cite{f:lc}, we will write the resulting formulas for the opposite
choice of the Borel subalgebra, as compared to the homomorphism \eqref{yhchcu}.
To this end, in types $B$, $C$ and $D$
we consider the automorphism $\si$ of the Lie algebra $t^{-1}\g_N[t^{-1}]$
defined on the generators by
\beql{autosi}
\si: F_{ij}[r]\mapsto -F_{ji}[r].
\eeq
We get the commutative diagram
\beql{cdhchne}
\begin{CD}
\U\big(t^{-1}\g_N[t^{-1}]\big)^{\h} @>>> \CC\big[\mu_i[r]\big]\\
@V\si VV     @VV\si V\\
\U\big(t^{-1}\g_N[t^{-1}]\big)^{\h} @>\chi>> \CC\big[\mu_i[r]\big],
\end{CD}
\eeq
where $i$ ranges over the set $\{1,\dots,n\}$ while $r< 0$.
The top and bottom horizontal arrows indicate the versions
of the Harish-Chandra homomorphism defined as in
\eqref{yhchcu}, where the left ideal $I$ is now
generated by
all elements $F_{ij}[r]$ with the conditions
$1\leqslant i<j\leqslant N$ and $r< 0$ for the top arrow, and by
all elements $F_{ij}[r]$ with the conditions
$N\geqslant i>j\geqslant 1$ and $r< 0$ for the bottom arrow
(which we denote by $\chi$).
The second vertical arrow indicates the isomorphism which takes
$\mu_i[r]$ to $-\mu_i[r]$.

Note that an automorphism analogous to \eqref{autosi} can be used in the case
of the Lie algebra $\gl_N$ to get the corresponding description
of the homomorphism $\chi$ and to
derive the formulas \eqref{hcha} and \eqref{hchh}.
However, these formulas follow easily
from the observation that $\tau+E[-1]$ is a Manin matrix
by the same argument as in the proof of Proposition~\ref{prop:hchan}.

To state the result in types $B$, $C$ and $D$, introduce the formal series
\beql{srene}
\ga_m(\om)\ts\tr\ts S^{(m)} \big(\di_u+F_1(u)_+\big)\dots \big(\di_u+F_m(u)_+\big),
\eeq
where we use notation
\eqref{galm} with $\om=N$ and $\om=-N$ in the orthogonal and symplectic case,
respectively, and
\ben
F(u)_+=\sum_{i,j=1}^N e_{ij}\ot F_{ij}(u)_+\in\End\CC^{N}\ot
\U\big(t^{-1}\g_N[t^{-1}]\big)[[u]].
\een
We will assume that in the symplectic case
the values of $m$ in \eqref{srene} are restricted to
$1\leqslant m\leqslant 2n+1$; see \cite[Sec.~3.3 and Sec.~4.1]{m:ff}.
The trace
is taken over all $m$ copies $\End\CC^N$
in the algebra
\beql{tenprcuu}
\underbrace{\End\CC^{N}\ot\dots\ot\End\CC^{N}}_m{}\ot
\U\big(t^{-1}\g_N[t^{-1}]\big)[[u,\di_u]]
\eeq
and we use matrix notation as in \eqref{matnot}. We set
\ben
\mu_i(u)_+=\sum_{r=0}^{\infty} \mu_i[-r-1]\ts u^r, \qquad i=1,\dots,n.
\een

\bpr\label{prop:hchnp}
The image of the series \eqref{srene} under
the homomorphism $\chi$
is given by the formula{\tss\rm:}
\ben
\bal
&\text{type $B_n${\rm:}\qquad\qquad} h_m\big(\di_u+\mu_{1}(u)_+,\dots,
\di_u+\mu_{n}(u)_+,\di_u-\mu_{n}(u)_+,\dots
\di_u-\mu_{1}(u)_+\big),\\[1.5em]
&\text{type $D_n${\rm:}\qquad\qquad}{\textstyle \frac{1}{2}}\ts h_m\big(\di_u+\mu_{1}(u)_+,\dots,
\di_u+\mu_{n-1}(u)_+,\di_u-\mu_{n}(u)_+,\dots
\di_u-\mu_{1}(u)_+\big)\\[0.5em]
{}&{\qquad\qquad\qquad\quad}\quad+{\textstyle \frac{1}{2}}\ts h_m\big(\di_u+\mu_{1}(u)_+,\dots,
\di_u+\mu_{n}(u)_+,\di_u-\mu_{n-1}(u)_+,\dots
\di_u-\mu_{1}(u)_+\big),\\[1.5em]
&\text{type $C_n${\rm:}\qquad\qquad}e_m\big(\di_u+\mu_{1}(u)_+,\dots,
\di_u+\mu_{n}(u)_+,\di_u,\di_u-\mu_{n}(u)_+,\dots
\di_u-\mu_{1}(u)_+\big).
\eal
\een
\epr

\bpf
We start with the orthogonal case $\g_N=\oa_N$.
The argument in the beginning of this section shows that
the image of the series
\ben
\ga_m(N)\ts\tr\ts S^{(m)} \big(\di_u-F_1(u)_+\big)\dots \big(\di_u-F_m(u)_+\big)
\een
under the homomorphism given by the top horizontal arrow in \eqref{cdhchne}
is found by Theorems~\ref{thm:hchbn} and \ref{thm:hchdn}, where
$\mu_i(u)$ should be respectively replaced by $-\mu_{i}(u)_+$
for $i=1,\dots,n$. Therefore, using the diagram \eqref{cdhchne}
we find that the image of the series
\beql{transpo}
\ga_m(N)\ts\tr\ts S^{(m)} \big(\di_u+F^{\tss t}_1(u)_+\big)
\dots \big(\di_u+F^{\tss t}_m(u)_+\big)
\eeq
under the homomorphism $\chi$ is given by the
respective $B_n$ and $D_n$ type formulas in the proposition,
where we set $F^{\tss t}(u)_+=\sum_{i,j}e_{ij}\ot F_{ji}(u)_+$.
It remains to observe that the series \eqref{transpo} coincides with \eqref{srene}.
This follows by applying the simultaneous transpositions
$e_{ij}\mapsto e_{ji}$ to all $m$ copies of $\End\CC^N$ and taking
into account the fact that $S^{(m)}$ stays invariant.

In the symplectic case, we suppose first that $m\leqslant n$.
Starting with the Harish-Chandra image provided by
Theorem~\ref{thm:hchcn} and applying
the same argument as in the orthogonal case,
we conclude that the image of the series
\beql{addser}
\ga_m(-2n)\ts\tr\ts S^{(m)} \big(\di_u-F_1(u)_+\big)\dots \big(\di_u-F_m(u)_+\big)
\eeq
under the homomorphism $\chi$ agrees
with the $C_n$ type formula given by the statement of the proposition. One more step
here is to observe that this series coincides with \eqref{srene}.
Indeed, this follows by applying
the simultaneous transpositions
$e_{ij}\mapsto \ve_i\tss\ve_j\tss e_{j'i'}$ to all $m$ copies of $\End\CC^N$.
On the one hand, this transformation does not affect the trace
of any element of \eqref{tenprcuu}, while on
the other hand, each factor $\di_u-F_i(u)_+$ is taken to $\di_u+F_i(u)_+$ and
the operator $S^{(m)}$ stays invariant.

Finally, extending the argument of \cite[Sec.~3.3]{m:ff} to the case $m=2n+1$
and using the results of \cite[Sec.~5]{m:ff}, we find that for all values
$1\leqslant m\leqslant 2n+1$
the coefficients $\Phi_{m\tss a}^{(s)}$ in the expansion
\ben
\ga_m(-2n)\ts\tr\ts S^{(m)} \big(\di_u+F_1(u)_+\big)\dots \big(\di_u+F_m(u)_+\big)
=\sum_{a=0}^m\sum_{s=0}^{\infty} \Phi_{m\tss a}^{(s)}\ts u^s\tss \di_u^{\ts a}
\een
belong to the Feigin--Frenkel center $\z(\wh\spa_{2n})$. The image of
the element $\Phi_{m\tss a}^{(s)}$ under the isomorphism
\eqref{hchiaff} is a polynomial in the generators $T^{\tss r}\Ec_{2k}$
of the classical $\Wc$-algebra $\Wc(\oa_{2n+1})$, where $k=1,\dots,n$ and $r\geqslant 0$;
see \eqref{bmglnb}. For a fixed value of $m$ and varying values of $n$
the coefficients of the polynomial are rational functions in $n$. Therefore, they
are uniquely determined by infinitely many values of $n\geqslant m$. This allows us
to extend the range of $n$ to all values $n\geqslant (m-1)/2$ for which
the expression \eqref{srene} is defined.
\epf

\bco\label{cor:extva}
Theorem~\ref{thm:hchcn} holds for all values $1\leqslant m\leqslant 2n+1$.
\eco

\bpf
This follows by reversing the argument used in the proof of Proposition~\ref{prop:hchnp}.
\epf

With the exception of the formula \eqref{pfaim} for the image
of the element $\phi'_n$ in type $D_n$, all
statements of the Main Theorem now follow from Proposition~\ref{prop:hchnp}.
It suffices to note that the coefficients of the polynomial \eqref{deftr}
and the differential operator \eqref{srene} are related via the vertex
algebra structure on the vacuum module $V_{-h^{\vee}}(\g_N)$. In particular,
the evaluation of the coefficients of the differential operator \eqref{srene}
at $u=0$ reproduces the corresponding coefficients of the polynomial \eqref{deftr}.
This implies the desired formulas for the Harish-Chandra images
in the Main Theorem; see e.g. \cite[Ch.~2]{f:lc} for the relevant properties
of vertex algebras.

Now consider the element $\Ec^{\tss\prime}_n$ of the algebra $\Wc(\oa_{2n})$
defined in \eqref{pfaffim} and which coincides with
the element \eqref{pfaim}.
To prove that the Harish-Chandra image of
the element $\phi'_n$ introduced by \eqref{genpf} equals
$\Ec^{\tss\prime}_n$,
use the automorphism of the Lie algebra $t^{-1}\oa_{2n}[t^{-1}]$
defined on the generators by
\beql{extau}
F^{}_{k\tss l}[r]\mapsto F_{\tilde k\tss \tilde l}[r],
\eeq
where $k\mapsto \tilde k$ is the involution on the set $\{1,\dots,2n\}$
such that $n\mapsto n'$, $n'\mapsto n$
and $k\mapsto k$ for all $k\ne n,n'$.
Note that $\phi'_n\mapsto -\phi'_n$ under the automorphism \eqref{extau}.
Similarly, $\Ec^{\tss\prime}_n\mapsto -\Ec^{\tss\prime}_n$ with respect to
the automorphism of $t^{-1}\h_{2n}[t^{-1}]$
induced by \eqref{extau}.

As a corollary
of the Main Theorem and the results of \cite{m:ff} we obtain
from the isomorphism \eqref{hchiaff} that
the elements
\begin{multline}
\Fc_m={\textstyle \frac{1}{2}}\ts h_m\big(T+\mu_{1}[-1],\dots,
T+\mu_{n-1}[-1],T-\mu_{n}[-1],\dots
T-\mu_{1}[-1]\big)\\[0.7em]
{}+{\textstyle \frac{1}{2}}\ts h_m\big(T+\mu_{1}[-1],\dots,
T+\mu_{n}[-1],T-\mu_{n-1}[-1],\dots
T-\mu_{1}[-1]\big),
\non
\end{multline}
with $m=2,4,\dots,2n-2$ together with
$\Ec^{\tss\prime}_n$
form a complete set of generators of $\Wc(\oa_{2n})$
(this fact does not rely on the calculation of the image of the Pfaffian).
Observe that all elements $T^{\tss r}\Fc_{2k}$
with $k=1,\dots,n-1$ and $r\geqslant 0$ are stable
under the automorphism \eqref{extau}. Since the Harish-Chandra image $\chi(\phi'_n)$
is a unique polynomial in the generators of $\Wc(\oa_{2n})$ and its degree
with respect to the variables $\mu_1[-1],\dots,\mu_n[-1]$ does not exceed $n$,
we can conclude that $\chi(\phi'_n)$ must be proportional to $\Ec^{\tss\prime}_n$.
The coefficient of the product $\mu_1[-1]\dots \mu_n[-1]$ in each of these two
polynomials is equal to $1$ thus proving that $\chi(\phi'_n)=\Ec^{\tss\prime}_n$.
This completes the proof of the Main Theorem.

The properties of vertex algebras mentioned above
and the relation $\chi(\phi'_n)=\Ec^{\tss\prime}_n$ imply
the respective formulas for the Harish-Chandra images
of the Pfaffians $\Pf\ts\wt F(u)_+$ and $\Pf\ts\wt F(u)$
defined by \eqref{genpf} with the matrix $\wt F[-1]$ replaced
by the skew-symmetric matrices $\wt F(u)_+=[F_{ij'}(u)_+]$
and $\wt F(u)=[F_{ij'}(u)]$, respectively.

\bco\label{cor:pfaim}
The Harish-Chandra images
of the Pfaffians are found by
\ben
\bal
\chi:\Pf\ts\wt F(u)_+&\mapsto
\big(\mu_{1}(u)_+-\di_u\big)\dots \big(\mu_{n}(u)_+-\di_u\big)\ts 1,\\[0.5em]
\Pf\ts\wt F(u)&\mapsto
\big(\mu_{1}(u)-\di_u\big)\dots \big(\mu_{n}(u)-\di_u\big)\ts 1,
\eal
\een
where the second map is defined in \eqref{yhchcu}.
\eco

\bpf
The first relation follows by the application of the state-field correspondence map
to the Segal--Sugawara vector \eqref{genpf} and using its
Harish-Chandra image \eqref{pfaim}. To get the second relation,
apply the automorphism \eqref{autosi} to the first relation to
calculate the image of $\Pf\ts\wt F(u)_+$ with respect to the homomorphism defined
by the top arrow in \eqref{cdhchne},
\ben
\Pf\ts\wt F(u)_+\mapsto
\big(\mu_{1}(u)_++\di_u\big)\dots \big(\mu_{n}(u)_++\di_u\big)\ts 1.
\een
Now replace $\wt F(u)_+$ with $-\wt F(u)$ and replace $\mu_{i}(u)_+$ with $-\mu_{i}(u)$
for $i=1,\dots,n$.
\epf

The isomorphism \eqref{hchiaff} and the Main Theorem provide
complete sets of generators of the classical $\Wc$-algebras.
In types $B$ and $C$ they coincide with those introduced in
Sec.~\ref{subsec:wabcd}, but different in type $D$, as pointed out
in the above argument.

\bco\label{cor:agen}
The elements $\Fc_2,\Fc_4,\dots,\Fc_{2n-2},\Ec^{\tss\prime}_n$
form a complete set of generators of $\Wc(\oa_{2n})$.
\qed
\eco

To complete this section, we point out that the application of the state-field
correspondence map to the coefficients of the polynomial \eqref{deftr}
and to the additional element \eqref{genpf} in type $D_n$ yields
Sugawara operators associated with $\wh\g_N$.
They act as scalars in the Wakimoto modules at the critical level.
The eigenvalues are found from the respective formulas of Proposition~\ref{prop:hchnp}
and Corollary~\ref{cor:pfaim} as follows from
the general theory of Wakimoto modules and their connection with the
classical $\Wc$-algebras; see \cite[Ch.~8]{f:lc}.

\section{Casimir elements for $\g_N$}
\label{sec:ce}
\setcounter{equation}{0}

We apply the theorems of Sec.~\ref{sec:ic} to
calculate the Harish-Chandra images of certain Casimir elements
for the orthogonal and symplectic Lie algebras previously considered in \cite{imr:ce}.
Our formulas for the Harish-Chandra images are equivalent to
those in \cite{imr:ce}, but take a different form.
We will work with the isomorphism \eqref{hchicl}, where the Cartan subalgebra
$\h$ of the Lie algebra $\g=\g_N$ is defined in the beginning of Sec.~\ref{sec:cy}
and the subalgebra $\n_+$ is spanned by the elements $F_{ij}$ with
$1\leqslant i<j\leqslant N$. We will use the notation $\mu_i=F_{ii}$ for
$i=1,\dots,N$ so that $\mu_i+\mu_{i'}=0$ for all $i$.

Consider the evaluation homomorphism
\ben
\ev: \U(\g_N[t])\to \U(\g_N),\qquad F_{ij}(u)\mapsto F_{ij}\tss u^{-1},
\een
so that $F_{ij}[0]\mapsto F_{ij}$ and $F_{ij}[r]\mapsto 0$ for $r\geqslant 1$.
The image of the series $\mu_i(u)$ then coincides with $\mu_i\ts u^{-1}$.
Applying the evaluation homomorphism to the series involved
in Theorems~\ref{thm:hchbn},
\ref{thm:hchdn} and \ref{thm:hchcn} we get the corresponding
Harish-Chandra images of the elements of the center
of the universal enveloping algebra $\U(\g_N)$. The formulas are obtained
by replacing $F_{ij}(u)$ with $F_{ij}\tss u^{-1}$ and $\mu_i(u)$ with $\mu_i\ts u^{-1}$.
Multiply the resulting formulas by $u^m$ from the left.
In the case $\g_N=\oa_{2n+1}$ use the relation
\beql{multum}
u^m\ts (\di_u+F_1\tss u^{-1})\dots (\di_u+F_m\tss u^{-1})
=(u\tss\di_u+F_1-m+1)\dots (u\tss\di_u+F_m)
\eeq
to conclude that
the Harish-Chandra image of the polynomial
\beql{casgenvb}
\ga_m(N)\ts\tr\ts S^{(m)} (F_1+v-m+1)\dots (F_m+v)
\eeq
with $v=u\tss\di_u$ is found by
\ben
\sum_{1\leqslant i_1\leqslant\dots\leqslant i_m\leqslant 1'}
(\mu_{i_1}+v-m+1)\dots (\mu_{i_m}+v),
\een
summed over the multisets $\{i_1,\dots,i_{m}\}$
with entries from $\{1,\dots,n,n',\dots,1'\}$.
By the arguments of \cite{imr:ce},
the Harish-Chandra image of the polynomial \eqref{casgenvb}
is essentially determined by those for the
even values $m=2k$
and a particular value of $v$.

\bco\label{cor:casb}
For $\g_N=\oa_{2n+1}$ the image of the Casimir element
\ben
\ga_{2k}(N)\ts\tr\ts S^{(2k)}\ts(F_1-k)\dots (F_{2k}+k-1)
\een
under the Harish-Chandra isomorphism
is given by
\beql{immub}
\sum_{1\leqslant i_1\leqslant\dots\leqslant i_{2k}\leqslant 1'}
(\mu_{i_1}-k)\dots (\mu_{i_{2k}}+k-1),
\eeq
summed over the multisets $\{i_1,\dots,i_{2k}\}$
with entries from $\{1,\dots,n,n',\dots,1'\}$.
Moreover, the element \eqref{immub} coincides with
the factorial complete symmetric function
\beql{imfacb}
\sum_{1\leqslant j_1\leqslant\dots\leqslant j_k\leqslant n}
\big(l^2_{j_1}-(j_1-1/2)^2\big)
\dots \big(l^2_{j_k}-(j_k+k-3/2)^2\big),
\eeq
where $l_i=\mu_i+n-i+1/2$ for $i=1,\dots,n$.
\eco

\bpf
The coincidence
of the elements \eqref{immub} and \eqref{imfacb} is verified
by using the characterization theorem for the factorial symmetric functions \cite{o:qi};
see also \cite{imr:ce}.
Namely, both elements are symmetric polynomials in $l_1^2,\dots,l_n^2$
of degree $k$, and their top degree components are both equal to the complete symmetric
polynomial $h_k(l_1^2,\dots,l_n^2)$. It remains to verify that each of
the elements \eqref{immub} and \eqref{imfacb} vanishes when $(\mu_1,\dots,\mu_n)$
is specialized to a partition with $\mu_1+\dots+\mu_n<k$ which
is straightforward.
\epf

Similarly, if $\g_N=\oa_{2n}$ use the same relation
\eqref{multum}
to conclude from Theorem~\ref{thm:hchdn} that
the Harish-Chandra image of the polynomial
\ben
2\ts\ga_m(N)\ts\tr\ts S^{(m)} (F_1+v-m+1)\dots (F_m+v)
\een
is found by
\ben
\sum_{\underset{\scriptstyle i_s\ne n}
{1\leqslant i_1\leqslant\dots\leqslant i_m\leqslant 2n}}
(\mu_{i_1}+v-m+1)\dots (\mu_{i_m}+v)
+
\sum_{\underset{\scriptstyle i_s\ne n'}
{1\leqslant i_1\leqslant\dots\leqslant i_m\leqslant 2n}}
(\mu_{i_1}+v-m+1)\dots (\mu_{i_m}+v),
\een
where the summation indices in the first sum do not include $n$
and the summation indices in the second sum do not include $n'$.

\bco\label{cor:casd}
For $\g_N=\oa_{2n}$ the image of the Casimir element
\ben
\ga_{2k}(N)\ts\tr\ts S^{(2k)}\ts(F_1-k)\dots (F_{2k}+k-1)
\een
under the Harish-Chandra isomorphism
is given by
\ben
{\textstyle\frac12}\sum_{\underset{\scriptstyle i_s\ne n}
{1\leqslant i_1\leqslant\dots\leqslant i_{2k}\leqslant 2n}}
(\mu_{i_1}-k)\dots (\mu_{i_{2k}}+k-1)
+
{\textstyle\frac12}\sum_{\underset{\scriptstyle i_s\ne n'}
{1\leqslant i_1\leqslant\dots\leqslant i_{2k}\leqslant 2n}}
(\mu_{i_1}-k)\dots (\mu_{i_{2k}}+k-1).
\een
Moreover, this element coincides with
the factorial complete symmetric function
\ben
\sum_{1\leqslant j_1\leqslant\dots\leqslant j_k\leqslant n}
\big(l^2_{j_1}-(j_1-1)^2\big)
\dots \big(l^2_{j_k}-(j_k+k-2)^2\big),
\een
where $l_i=\mu_i+n-i$ for $i=1,\dots,n$.
\eco

\bpf
The coincidence of the two expressions for the Harish-Chandra image
is verified in the same way as for the case of $\oa_{2n+1}$
outlined above.
\epf

Now suppose that $\g_N=\spa_{2n}$ and use the relation
\ben
u^m\ts (-\di_u+F_1\tss u^{-1})\dots (-\di_u+F_m\tss u^{-1})
=(-u\tss\di_u+F_1+m-1)\dots (-u\tss\di_u+F_m)
\een
to conclude from Theorem~\ref{thm:hchcn} and Corollary~\ref{cor:extva} that
the Harish-Chandra image of the polynomial
\ben
\ga_m(-2n)\ts\tr\ts S^{(m)} (F_1+v+m-1)\dots (F_m+v)
\een
with $v=-u\tss\di_u$ is found by
\ben
\sum_{1\leqslant i_1<\dots< i_m\leqslant\tss 1'}
(\mu_{i_1}+v+m-1)\dots (\mu_{i_m}+v),
\een
summed over the subsets $\{i_1,\dots,i_m\}$ of the set $\{1,\dots,n,0,n',\dots,1'\}$
with the ordering $1<\dots<n<0<n'<\dots<1'$, where $\mu_0:=0$.
Taking $m=2k$ and $v=-k+1$
we get the following.

\bco\label{cor:casc}
For $\g_N=\spa_{2n}$ the image of the Casimir element
\ben
\ga_{2k}(-2n)\ts\tr\ts S^{(2k)}\ts(F_1+k)\dots (F_{2k}-k+1)
\een
under the Harish-Chandra isomorphism
is given by
\beql{immuc}
\sum_{1\leqslant i_1<\dots< i_{2k}\leqslant 1'}
(\mu_{i_1}+k)\dots (\mu_{i_{2k}}-k+1),
\eeq
summed over the subsets $\{i_1,\dots,i_{2k}\}
\subset\{1,\dots,n,0,n',\dots,1'\}$.
Moreover, the element \eqref{immuc} coincides with
the factorial elementary symmetric function
\beql{imfacc}
(-1)^k\sum_{1\leqslant j_1<\dots< j_k\leqslant n}
\big(l^2_{j_1}-j_1^2\big)
\dots \big(l^2_{j_k}-(j_k-k+1)^2\big),
\eeq
where $l_i=\mu_i+n-i+1$ for $i=1,\dots,n$.
\eco

\bpf
To verify that
the elements \eqref{immuc} and \eqref{imfacc} coincide, use again
the characterization theorem for the factorial symmetric functions \cite{o:qi};
see also \cite{imr:ce}.
Both elements are symmetric polynomials in $l_1^2,\dots,l_n^2$
of degree $k$, and their top degree components are both equal to the elementary symmetric
polynomial $(-1)^k\tss e_k(l_1^2,\dots,l_n^2)$. Furthermore, it is easily seen that each of
the elements \eqref{immuc} and \eqref{imfacc} vanishes when $(\mu_1,\dots,\mu_n)$
is specialized to a partition with $\mu_1+\dots+\mu_n<k$.
\epf

\end{document}